

\input epsf.tex

\def\2{{1\over 2}}

\def\d{\delta}
\def\a{\alpha}
\def\b{\beta}
\def\g{\gamma}

\def\s{\sigma}
\def\e{\epsilon}
\def\l{\lambda}
\def\o{\omega}

\def\fun#1#2#3{#1\colon #2\rightarrow #3}

\def\frac#1#2{{{#1} \over {#2}}}

\def\sqr{\sqrt}
\def\st{\;\colon\;}
\def\tends{\rightarrow}

\def\dr{ {\rm d} }

\def\R{{\bf R}}
\def\N{{\bf N}}

\def\thm#1{\vskip 1 pc\noindent{\bf Theorem #1.\quad}\sl}
\def\lem#1{\vskip 1 pc\noindent{\bf Lemma #1.\quad}\sl}
\def\prop#1{\vskip 1 pc\noindent{\bf Proposition #1.\quad}\sl}
\def\cor#1{\vskip 1 pc\noindent{\bf Corollary #1.\quad}\sl}

\def\proof{\rm\vskip 1 pc\noindent{\bf Proof.\quad}}
\def\fin{\par\hfill $\backslash\backslash\backslash$\vskip 1 pc}
\def\txt#1{\quad\hbox{#1}\quad}

\def\L{{\cal L}}

\def\s{\sigma}

\def\o{\omega}

\def\2{\frac{1}{2}}
\def\inn#1#2{{\langle #1 ,#2\rangle}}

\def\part{{\partial_{x}}}

\def\tr{{{}^{t}}}

\def\cc{{\cal C}}
\def\ec{{\cal E}}

\def\fc{{\cal F}}

\def\dc{{\cal D}}
\def\oc{{\cal O}}
\def\pc{{\cal P}}

\def\mc{{\cal M}}



\baselineskip= 17.2pt plus 0.6pt
\font\titlefont=cmr17
\centerline{\titlefont Another point of view}
\vskip 1 pc
\centerline{\titlefont on Kusuoka's measure}
\vskip 4pc
\font\titlefont=cmr12
\centerline{         \titlefont {Ugo Bessi}\footnote*{{\rm 
Dipartimento di Matematica, Universit\`a\ Roma Tre, Largo S. 
Leonardo Murialdo, 00146 Roma, Italy.}}   }{}\footnote{}{
{{\tt email:} {\tt bessi@matrm3.mat.uniroma3.it} Work partially supported by the PRIN2009 grant "Critical Point Theory and Perturbative Methods for Nonlinear Differential Equations}} 
\vskip 0.5 pc
 
\par
\vskip 2pc
\centerline{\bf Abstract}

Kusuoka's measure on fractals is a Gibbs measure of a very special kind, because its potential is discontinuous, while the standard theory of Gibbs measures requires continuous (actuallly, H\"older) potentials. In this paper, we shall see that for many fractals it is possible to build a class of matrix-valued Gibbs measures completely within the scope of the standard theory; there are naturally some minor modifications, but they are only due to the fact that we are dealing with matrix-valued functions and measures. We shall use these matrix-valued Gibbs measures to build self-similar Dirichlet forms on fractals. Moreover, we shall see that Kusuoka's measure can be recovered in a simple way from the matrix-valued Gibbs measure.

\vskip 2 pc
\centerline{\bf  Introduction}
\vskip 1 pc

First of all, let us briefly explain what we mean by a fractal $G$ on $\R^d$; our definition is less general than the one in [12]. We consider $n$ contractions
$$\psi_1,\dots,\psi_n\in C^{1,\nu_0}(\R^d,\R^d)    \eqno (1)  $$
with $\nu_0\in(0,1]$; it is standard ([9]) that there is a unique compact set $G\subset\R^d$ such that 
$$G=\bigcup_{i=1}^n\psi_i(G)  .  $$
We shall also suppose that the maps $\psi_i$ are the "branches of the inverse" of a Borel map $\fun{F}{G}{G}$. Since the maps $\{ \psi_i \}_{i=1}^n$ are contractions, their inverse $F$ is expanding; in Dynamical Systems, expanding maps have been studied extensively (see for instance [13] or [19]); as we shall see, many results on expanding maps carry over to Dirichlet forms on fractals (see [9] or [17] for an introduction to this theory). 

Applying ergodic theory to the study of Dirichlet forms on fractals is not new: in section 1.4 of [7] the idea of applying the Ruelle operator to the study of Kusuoka's measure is attributed to Strichartz. Actually, in this paper we try to understand the results of [7] and [12] by looking at them from a slightly different perspective. 

In order to explain the connection between expanding maps and Dirichlet forms, we begin recalling the scalar Gibbs measure; though it will play no r\^ole in our paper, it will guide us in the construction of the matrix-valued one. 

For $\nu_0\in (0,1]$ we take 
$V\in C^{\nu_0}(G,\R)$ and define the scalar Ruelle operator as 
$$\fun{\L_{sc}}{C(G,\R)}{C(G,\R)},\qquad
(\L_{sc} v)(x)=
\sum_{i=1}^ne^{V\circ\psi_i(x)}v\circ\psi_i(x)   .   \eqno (2)$$

Using the Perron-Frobenius theorem one can prove ([18]) that there is $\b>0$ and a continuous function $h>0$ such that $\L_{sc} h=\b h$. Since $\L_{sc}$ is a continuous operator from the space of continuous functions into itself, its adjoint $\L_{sc}^\ast$ brings the space of Borel measures into itself; it can be shown that there is a measure $\mu$, called the Gibbs measure, such that $h\mu$ is probability and $\L^\ast_{sc}\mu=\b\mu$. One of the properties of $\mu$ is the following: for all $u,v\in C(G,\R)$ we have that 
$$\frac{1}{\b}\int_G u(\L_{sc} v)\dr\mu=
\int_G(u\circ F)v\dr\mu   .    \eqno (3)$$

We point out a few consequences of (3). First of all, we can write the adjoint $\L_{sc}^\ast$ explicitly, at least for measures absolutely continuous with respect to $\mu$.   
$$\frac{1}{\b}\L_{sc}^\ast(u\mu)=(u\circ F)\mu  .  $$
Moreover, (3) tells us that $\frac{1}{\b}\L_{sc}v$ is the density of $F_\sharp(v\mu)$, where 
$F_\sharp\nu$ denotes the push-forward of the measure $\nu$ by $F$; the push forward is defined by 
$$\int_Gf\dr(F_\sharp\nu)=\int_Gf\circ F\dr\nu  $$
for all $f\in C(G,\R)$. Since $\L_{sc}h=\b h$, this implies that $F_\sharp(h\mu)=h\mu$, i. e. that $h\mu$ is an invariant measure. 

Moreover, we know how $\mu[\psi_{i_0}\circ\dots\psi_{i_l}(G)]$ scales as $l\tends+\infty$; namely, there is a constant $D_1>0$, independent of the sequence $i_0i_1\dots$, such that 
$$\frac{1}{D_1}\le
\frac{
\mu[\psi_{i_0}\circ\dots\psi_{i_l}(G)]
}{
\b^{-l}\cdot\exp(V(x)+V(F(x))+\dots+V(F^{l-1}(x)))
}   \le D_1   $$
for all $x\in \psi_{i_0}\circ\dots\psi_{i_l}(G)$. 

Kusuoka's measure $\kappa$ is defined by a similar scaling property; when the maps 
$\psi_i$ of (1) are affine, i. e. $D\psi_i$ is a constant matrix, we have  
$$\kappa[(\psi_{x_0}\circ\dots\psi_{x_l})(G)]=
\frac{1}{\b^l}\cdot{\rm tr}[\hat Q (D\psi_{x_0}\cdot\dots\cdot D\psi_{x_l})Q
\tr(D\psi_{x_0}\cdot\dots\cdot D\psi_{x_l})]$$
where $Q$, $\hat Q$ are two suitable symmetric $d\times d$ matrices, $\b>0$, $\tr Q$ denotes the transpose of $Q$ and ${\rm tr}$ is the trace.  As explained in [7] (see [3] for the proof and further details), $\kappa$ is a Gibbs measure, but for a discontinuous potential $V$, and the standard theory does not apply to it. 

Let us come to the Dirichlet form. On many fractals, the space $L^2(S,\kappa)$ admits a "heat semigroup" $P_s$ which is induced, as in $\R^d$, by a Brownian motion ([1], [2], [6], [11]). The semigroup $P_s$ has a generator $A$ 
$$Af=\lim_{h\searrow 0}\frac{P_hf-f}{h}$$
which induces a Dirichlet form on $L^2(G,\kappa)$; this form is defined on a dense subspace $\dc(\ec)\subset L^2(S,m)$ by 
$$\ec(f,g)=-\int_G (Af)\cdot g\dr m  \txt{if}g\in\dc(\ec)
\txt{and}f\in\dc(A)   .   $$
Conversely ([5]), given a Dirichlet form $\ec$ it is relatively easy to check whether it is induced by a Brownian motion; since Dirichlet forms are easier to study than the Brownian motion itself, they immediately attracted attention ([14], [15]; a counterexample to the existence of Dirichlet forms is in [16]). 

This brings us to the matrix-valued Gibbs measure: under suitable hypotheses on the fractal $G\subset\R^d$, the "natural" Dirichlet form $\ec$ on $L^2(G,\kappa)$ can be written in the following way: if $u,v\in C^1(\R^d,\R)$, then 
$$\ec(u,v)=\int_G (T_{x}\nabla u(x),\nabla v(x))\dr \kappa(x)$$
where $(\cdot,\cdot)$ denotes the standard inner product in $\R^d$, $T_x$ is a Borel field of symmetric matrices (in many cases, projections) and $\kappa$ is Kusuoka's measure. 

The matrix-valued measure $\tau\colon=T_{x}\kappa$ appears in a natural way in the formula above; the aim of this paper is to show that $\tau$ is a Gibbs measure as well. 

More precisely, we denote by $M^d$ the space of symmetric $d\times d$ matrices; we can define a Ruelle operator
$$\fun{\L_G}{C(G,M^d)}{C(G,M^d)}$$
by
$$(\L_G A)(x)=
\sum_{i=1}^n\tr D\psi_i(x) A(\psi_i(x)) D\psi_i(x)  .  \eqno(4)$$
The dual space of $C(G,M^d)$ is the space $\mc(G,M^d)$ of $M^d$-valued measures on $G$, and the adjoint $\L^\ast_G$ of $\L_G$ brings $\mc(G,M^d)$ into itself. As we shall see, $\L_G$ and $\L^\ast_G$ have each a positive-definite eigenvector, which we shall call $Q_G$ and $\tau_G$ respectively. In lemma 4.8 below we shall prove the following version of (3): if $g\in C(G,\R)$ and $A\in C(G,M^d)$, defining the integral as in section 2 below we have 
$$\int_G\left(
g\cdot\left(
\frac{1}{\b}\L_G A
\right)   ,\dr\tau_G
\right)_{HS}=
\int_G(g\circ F\cdot A,\dr\tau_G)_{HS}  .  $$
The formula above implies that the scalar measure $(Q_G,\tau_G)_{HS}$ (again, see section 2 for the definition) is invariant. 

\thm{1} Let the maps $\psi_1,\dots,\psi_n$ and $F$ satisfy hypotheses (F1)-(F4) of section 1 below and (ND) with constant $b>0$ at the beginning of section 4. Let $G$ be the fractal associated with $\psi_1,\dots,\psi_n$ and let $M^d$ denote the space of $d\times d$ symmetric matrices. Let the operator $\L_G$ be defined as in (4)  and let   
$||D\psi_i||_{\nu_0}$ be smaller than a positive constant that depends on $b>0$ (which is always true if (ND) holds and the maps $\psi_i$ are affine); then the following holds.

\noindent 1) There are $Q_G\in C(G,M^d)$ and $\b>0$ such that 
$$(\L_G Q_G)(x)=\b Q_G(x)  \qquad\forall x\in G.  $$
The map $Q_G$ belongs to $C^{\nu_0}(G,M^d)$ and is unique up to multiplication by a constant; again up to multiplication by a constant, $Q_G(x)$ is positive-definite for all 
$x\in G$. 

\noindent 2) Let $\L^\ast_G$ denote the adjoint of $\L_G$; then, there is a Borel measure 
$\tau_G$ on $G$ which takes values in $M^d$ and such that
$$\L^\ast_G\tau_G=\b\tau_G  .  $$
The measure $\tau_G$ is unique up to multiplication by a constant; again up to multiplication by a constant, $\tau_G$ takes values in semi-positive definite matrices. 

\noindent 3) The measures $||\tau_G||$ and $\kappa_G\colon=(Q_G,\tau_G)_{HS}$ are mutually absolutely continuous. Moreover, $\kappa_G$ is ergodic for the map $F$. 

\noindent 4) We define the form $\fun{\ec}{C^1(\R^d)\times C^1(\R^d)}{\R}$ in the following way (the notation for the integral is in section 2 below): 
$$\ec(f,g)=
\int_G(\nabla f(x),\dr\tau_G(x)\nabla g(x))  .  $$
Then, $\ec$ is self-similar, i.e. 
$$\ec(f,g)=\frac{1}{\b}\sum_{i=1}^n\ec(f\circ\psi_i,g\circ\psi_i)     $$
for all $f,g\in C^1(\R^d)$. 

\noindent 5) Let the maps $\psi_i$ be affine. Then, the measure $\tau_G$ has the Gibbs property, i. e., for all $x_0,\dots,x_l\in(1,\dots,n)$ we have that 
$$\tau_G(\psi_{x_0}\circ\dots\circ\psi_{x_{l-1}}(G))=
\frac{1}{\b^l}\cdot (D(\psi_{x_0}\circ\dots\circ\psi_{x_{l-1}}))\cdot\tau_G(G)\cdot 
\tr (D(\psi_{x_0}\circ\dots\circ\psi_{x_{l-1}}))   .   $$
We haven't written explicitly the point where we calculate 
$D(\psi_{x_0}\circ\dots\circ\psi_{x_l})$ since these maps are affine and their derivative is constant.

\rm

\vskip 1pc

Note that in point 4) we do not assert that $\ec$ is closable: actually, we don't know any criteria for closability other than the ones in [9] and [12]. 

The paper is organised as follows. In section 1 we recall the notation and the basic facts about the Perron-Frobenius theorem, fractal sets and Dirichlet forms. In section 2 we define the convex cones to which we are going to apply the Perron-Frobenius theorem. In section 3 we define the Ruelle operator $\L_G$ on matrices and show that the fixed points of its adjoint $\L^\ast_G$ induce a self-similar quadratic form $\ec$ on $C^1(\R^d)$. In section 4, we apply the Perron-Frobenius theorem to find the maximal eigenvector of 
$\L_G$ and the matrix-valued Gibbs measure $\tau_G$. In section 5, we show that 
$\tau_G$ has the Gibbs property. 

\vskip 2pc

\centerline{\bf\S 1}

\centerline{\bf Preliminaries and notation}

\vskip 1pc

\noindent{\bf The Perron-Frobenius theorem.} We follow [19] (see also [4] for the original treatment). 

Let $X$ be a real vector space; we say that $\cc\subset X\setminus\{  0  \}$ is a cone if 
$$v\in\cc
\txt{and}
t>0
\txt{implies that}
tv\in\cc   .  $$
Let $\cc\subset X$ be a convex cone; we say that $w\in\bar\cc$ if there are $v\in\cc$ and 
$t_n\searrow 0$ such that $w+t_nv\in\cc$ for all $n\ge 1$. In what follows, we shall suppose that $\cc$ is a convex cone such that 
$$\bar\cc\cap(-\bar\cc)=\{ 0 \}  .  \eqno (1.1)$$
If $v_1,v_2\in\cc$, we define 
$$\a(v_1,v_2)=\sup\{
t>0\st v_2-tv_1\in\cc
\}   \eqno (1.2)$$
$$\frac{1}{\b(v_1,v_2)}=
\sup\{
t>0\st v_1-tv_2\in\cc
\}  \eqno (1.3)$$
and
$$\theta(v_1,v_2)=\log
\frac{\b(v_1,v_2)}{\a(v_1,v_2)}   .   \eqno (1.4)$$

Since $\theta(v,\l v)=0$ for all $\l>0$, we identify the points of a ray; namely, we say that $v_1\simeq v_2$ if $v_2=tv_1$ for some $t>0$; we shall denote by 
$\frac{\cc}{\simeq}$ the set of equivalence classes. 

We have that $\theta(v_1,v_2)\in[0,+\infty]$ for all $v_1,v_2\in C$; if $\theta$ never assumes the value $+\infty$, then $\theta$ is a distance on $\frac{\cc}{\simeq}$. 

The following proposition from [19] allows us to use the contraction principle. 

\prop{1.1} 1) Let $\fun{L}{X}{X}$ be a linear operator such that $L(\cc)\subset\cc$ and let us define
$$D=\sup\{
\theta(Lv_1,Lv_2)\st v_1,v_2\in\cc
\}   .    $$
Then, if $D<+\infty$, $L$ is a contraction on $(\frac{\cc}{\simeq},\theta)$, namely
$$\theta(Lv_1,Lv_2)\le (1-e^{-D})\theta(v_1,v_2)  \qquad
\forall v_1,v_2\in\cc  .  $$

\noindent 2) As a consequence of 1), if $D<+\infty$ and $(\frac{\cc}{\simeq},\theta)$ is a complete metric space, there is $(\l,v)\in(0,+\infty)\times\cc$, unique in 
$(0,+\infty)\times\frac{\cc}{\simeq}$, such that 
$$Lv=\l v  .  $$
Moreover, if $w\in\cc$, then 
$$\theta(L^nw,v)\le\theta(w,v)\frac{(1-e^{-D})^n}{e^{-D}}  .  \eqno (1.5)$$

\rm

\vskip 1pc

\noindent{\bf Fractal sets.} We make the following hypotheses on the fractal set.

\noindent {\bf(F1)} There is $\nu_0\in(0,1]$ and diffeomorphisms 
$$\psi_1,\dots,\psi_n\in C^{1,\nu_0}(\R^d,\R^d)   \eqno (1.6)$$
satisfying 
$$\eta\colon=\sup_{i\in(1,\dots,n)}Lip(\psi_i)<1  .  \eqno (1.7)$$
By theorem 1.1.7 of [9], this implies that 
there is a unique non empty compact set $G\subset\R^d$ such that
$$G=\bigcup_{i=1}^n\psi_i(G) . \eqno (1.8)$$
In the following, we shall always rescale the norm of $\R^d$ in such a way that 
$${\rm diam}(G)\le 1   .  \eqno (1.9)$$

If (F1) holds, then the dynamics of $F$ on $G$ can be coded. Indeed, we define 
$\Sigma$ as the space of sequences 
$$\Sigma=\{ 1,\dots,n \}^\N=
\{  \{ x_i \}_{i\ge 0}\st x_i\in(1,\dots,n),\quad\forall i\ge 0
\}    $$
with the product topology. This is a metric space; for instance, if $\g\in(0,1)$, we can define the metric 
$$d_\g(\{ x_i \}_{i\ge 0},\{ y_i \}_{i\ge 0})=\g^k$$
where 
$$k=\inf\{
i\ge 0\st x_i\not=y_i
\}   ,   $$
with the convention that the $\inf$ of the empty set is $+\infty$. 

We define the shift $\s$ as 
$$\fun{\s}{\Sigma}{\Sigma},\qquad
\fun{\sigma}{\{ x_0,x_1,x_2,\dots \}}{\{ x_1,x_2,x_3,\dots \}}  .   $$

If $x_0,\dots,x_l\in(1,\dots,n)$, we define the cylinder 
$$[x_0\dots x_l]=\{
\{ y_i \}_{i\ge 0}\st y_i=x_i\txt{for}i\in(1,\dots,l)
\}  .  $$
We also set 
$$\psi_{x_0\dots x_l}=\psi_{x_0}\circ\dots\circ\psi_{x_l}  $$
and
$$[x_0\dots x_l]_G=\psi_{x_0}\circ\psi_{x_{1}}\circ\dots\circ\psi_{x_l}(G)  .  \eqno (1.10)$$
If $x=(x_0x_1\dots)$ we set $(ix)=(ix_0x_1\dots)$. Now (1.10) implies that 
$$\psi_i([x_1\dots x_l]_G)=[i x_1\dots x_l]_G  .  \eqno (1.11)$$
Since the maps $\psi_i$ are continuous and $G$ is compact, the sets 
$[x_0\dots x_l]_G\subset G$ are compact. By (1.8) we have that $\psi_i(G)\subset G$ for 
$i\in(1,\dots,n)$; this implies that, for all 
$\{ x_i \}_{i\ge 0}\in\Sigma$
$$[x_0\dots x_{l-1} x_l]_G\subset [x_0\dots x_{l-1}]_G  .  $$
From (1.7), (1.9) and (1.10) we get that 
$${\rm diam} ([x_0\dots x_l]_G)\le\eta^l   .   \eqno (1.12)$$
Let $\{ x_i \}_{i\ge 0}\subset\Sigma$; by the last two formulas and the finite intersection property we have that 
$$\bigcap_{l\ge 1}[x_0\dots x_l]_G$$
is a single point, which we call $\Phi(\{ x_i \}_{i\ge 0})$; formula (1.12) implies in a standard way that the map $\fun{\Phi}{\Sigma}{G}$ is continuous. It is not hard to prove, using (1.8), that $\Phi$ is surjective. We shall call $\tilde d$ the distance on $G$ induced by the Euclidean distance on $\R^d$ and, from now on, in our choice of the metric on $\Sigma$ we take 
$\g\in(\eta,1)$; this implies by the definition of $d_\g$ and (1.12) that $\Phi$ is 1-Lipschitz.  

\noindent {\bf (F2)} If $i\not=j$, $\psi_i(G)\cap\psi_j(G)$ is a finite set. We set 
$$\fc\colon=\bigcup_{i\not =j}\psi_i(G)\cap\psi_j(G)  .  $$

\noindent {\bf (F3)} We ask that $G$ is post-critically finite, which means the following: if we set $A=\Phi^{-1}(\fc)$, then the set $\cup_{j\ge 1}\s^j(A)$ is finite.

\noindent {\bf (F4)} We ask that there are disjoint open sets $\oc_1,\dots,\oc_n\subset\R^d$ such that
$$G\cap\oc_i=\psi_i(G)\setminus\left(\bigcup_{i\not=j}\psi_i(G)\cap\psi_j(G)\right)\txt{for}i\in(1,\dots,n)  .  $$
We define a map $\fun{F}{\bigcup_{i=1}^n\oc_i}{\R^d}$ by
$$F(x)=\psi_i^{-1}(x)  \txt{if}x\in\oc_i  .  $$
If moreover we ask that $\oc_i\subset\psi_i^{-1}(\oc_i)$ (or, equivalently, that 
$\psi_i(\oc_i)\subset\oc_i$, since the maps $\psi_i$ are diffeos), this implies the first equality below.
$$F\circ\psi_i(x)=x\qquad\forall x\in\oc_i\subset \psi_i^{-1}(\oc_i)
\txt{and}
\psi_i\circ F(x)=x\qquad\forall x\in\oc_i  .  \eqno (1.13)$$
We call $a_i$ the unique fixed point of $\psi_i$; note that, by (1.8), $a_i\in G$. If $x\in\fc$, we define $F(x)=a_j$ for some arbitrary $a_j$.  This defines $F$ as a Borel map on all of $G$, which satisfies (1.13). 

We point out a consequence of (F4): the sets $\oc_i$ and $\psi_{j}(G)$ do not intersect unless $i=j$. By the definition of the coding this implies that, if 
$z=\Phi(\{ x_j \}_{j\ge 0})\in\oc_i$, then $x_0=i$. Since $\psi_{x_0}$ is a diffeo, also the sets $\psi_{x_0}(\oc_j)$ and $\psi_{x_0}\circ\psi_l(G)$ do not intersect unless $l=j$. Thus, if in addition $z\in\bigcup_{j=1}^n\psi_{x_0}(\oc_j)$, $z$ belongs to a unique 
$\psi_{x_0}(\psi_{x_1}(G))$ and also $x_1$ is uniquely determined. Going on, we see that $z$ has a unique coding unless it belongs to the countable set 
$$\bigcup_{n\ge 1}\bigcup_{i_0,\dots,i_n}\psi_{i_0\dots i_n}(\fc)  .  $$

By (1.10) and the definition of $\Phi$ we easily get that 
$[x_0\dots x_l]\subset \Phi^{-1}([x_0\dots x_l]_G)$; since the set where $\Phi$ is not injective is countable, we have that  
$$\sharp(
\Phi^{-1} ([x_0\dots x_l]_G)\setminus [x_0\dots x_l]
)  \le\sharp\N  .  \eqno (1.14)$$
Note that, if $x=(x_0x_1\dots)$, then by the definition of $\Phi$ 
$$\Phi\circ\s(x)=\bigcap_{l\ge 1}[x_1\dots x_l]_G  .  $$
The definition of $F$ implies the first equality below; if we suppose that 
$\Phi(x)\in\oc_{x_0}$ and recall that $F=\psi_{x_0}^{-1}$ on $\oc_{x_0}$ we get the middle one while the last equality comes from the formula above. 
$$F\circ\Phi(x)=F\left(
\bigcap_{l\ge 1}[x_0\dots x_l]_G
\right)  =    
\bigcap_{l\ge 1}[x_1\dots x_l]_G=
\Phi\circ\s(x)  .  $$
In other words, the first equality below holds save when $\Phi(x)\in\fc$. The second equality below follows for all $x\in G$ from (1.11). 
$$\left\{\matrix{
\Phi\circ\s(x)=F\circ\Phi(x)\txt{save possibly when $\Phi(x)\in \fc$,}\cr
\Phi(i,x) =\psi_i(\Phi(x))   \qquad\forall x\in\Sigma,
\quad\forall i\in (1,\dots,n)   .   
}
\right.       \eqno (1.15)  $$
In other words, up to a change of coordinates, shifting the coding one place to the left is the same as applying $F$. Iterating the first one of (1.15) we get that, for all $l\ge 1$, 
$$\Phi\circ\s^l(x)=F^l\circ\Phi(x)\txt{save possibly for} 
x\in \bigcup_{j\ge 0}\s^{-j}(\Phi^{-1}(\fc))   . \eqno (1.16)$$
Note that the union on the right is a countable set, since $\Phi^{-1}(\fc)$ is finite by (F3). 

A particular case we have in mind is the harmonic Sierpinski gasket on $\R^2$ ([8], [10]). We set
$$T_1=\left(
\matrix{
\frac{3}{5},&0\cr
0,&\frac{1}{5}
}
\right)  ,  \quad
T_2=\left(
\matrix{
\frac{3}{10},&\frac{\sqr 3}{10}\cr
\frac{\sqr 3}{10},&\frac{1}{2}
}
\right)  ,\quad
T_3=\left(
\matrix{
\frac{3}{10},&-\frac{\sqr 3}{10}\cr
-\frac{\sqr 3}{10},&\frac{1}{2}
}
\right)   ,  $$
$$A=\left(
\matrix{
0\cr
0
}
\right),\qquad
B=\left(
\matrix{
1\cr
\frac{1}{\sqrt 3}
}
\right) ,\qquad
C=\left(
\matrix{
1\cr
-\frac{1}{\sqrt 3}
}
\right)  $$
and 
$$\psi_1(x)=T_1(x),\quad \psi_2(x)=
B   +T_2\left(
x-B
\right),\quad
\psi_3(x)=
C   +T_3\left(
x-C
\right)  .  $$
Referring to the figure below, $\psi_1$ brings the triangle $ABC$ into $Abc$; 
$\psi_2$ brings $ABC$ into $Bac$ and $\psi_3$ brings $ABC$ into $Cba$. We take 
$\oc_1$, $\oc_2$, $\oc_3$ as three disjoint open sets which contain, respectively, the triangle $Abc$ minus $b,c$, $Bca$ minus $c,a$ and $Cba$ minus $a,b$. 

We define the map $F$ as 
$$F(x)=\psi_i^{-1}(x)
\txt{if}
x\in\oc_i $$
and we extend it as in (F4) on $\{ a,b,c \}$.

It is easy to check that the fractal $G$ generated by $\psi_1,\psi_2,\psi_3$ satisfies hypotheses (F1)-(F4) above; it is easy to check that it also satisfies (ND) of section 4 below. 

\vskip 1pc

\noindent {\bf The invariant Dirichlet form.} Let $(X,\tilde d,\nu)$ be a metric measure space; we suppose for simplicity that $(X,\tilde d)$ is compact and $\nu$ is probability. 

A Dirichlet form is a symmetric bilinear form 
$$\fun{\ec}{\dc(\ec)\times\dc(\ec)}{\R}$$
defined on a dense set $\dc(\ec)\subset L^2(X,\nu)$ such that the two conditions below hold. 

\noindent (D1) $\dc(\ec)$ is closed under the graph norm; in other words, $\dc(\ec)$ is a Hilbert space for the norm
$$||u||_{\dc(\ec)}^2=||u||^2_{L^2(S,m)}+\ec(u,u)  .  \eqno (1.17)$$

\noindent (D2) $\ec$ is Markovian, i. e. 
$$\ec(\eta\circ f,\eta\circ f)\le\ec(f,f)$$
for all $f\in\dc(\ec)$ and all 1-Lipschitz maps $\fun{\eta}{\R}{\R}$ with $\eta(0)=0$.

We list some additional properties a Dirichlet form can have.  

\noindent (D3) $\ec$ is regular; this means that $\dc(\ec)\cap C(X,\R)$ is dense in 
$C(X,\R)$ for the uniform topology and in $\dc(\ec)$ for the graph norm (1.17). 

\noindent (D4) $\ec$ is strongly local, i. e. 
$$\ec(f,g)=0$$
whenever $f,g\in\dc(\ec)$ and $f$ is constant in a neighbourhood of the support of $g$. 

It can be proven ([5]) that, if $\ec$ satisfies (D1)-(D4), then it is the Dirichlet form of a Brownian motion. 

\noindent (D5) $\ec$ is self similar, i. e. there is $\b>0$ such that 
$$\ec(u,v)=\b\sum_{i=1}^n\ec(u\circ\psi_i,v\circ\psi_i)  $$
for all $u,v\in\dc(\ec)$. 

As we stated in the introduction, we shall be able to build a local and self-similar form on $C^1(\R^d,\R)$, but not to prove its closability.

\vskip 2pc
\centerline{\bf \S 2}
\centerline{\bf Function spaces and cones}
\vskip 1pc

We begin listing three equivalent ways to define the norm of a matrix in $M^d$. 

The first one is the $\sup$ norm
$$||A||=\sup\{
||Av||\st ||v||\le 1
\}   .  $$
We denote by $\tr A$ the adjoint of the matrix $A$ and by ${\rm tr}(A)$ its trace; on the space of all matrices we can define the inner product 
$$(A,B)_{HS}={\rm tr}(\tr AB)   $$
which induces the Hilbert-Schmidt norm
$$||A||_{HS}^2={\rm tr}(\tr A A)   .   $$
If $A$ is symmetric we can set 
$$|||A|||=\sup\{
|(Av,v)|\st ||v||\le 1
\}  $$
where we have denoted by $(\cdot,\cdot)$ the inner product of $\R^d$. 

Clearly, if $A$ is symmetric, $|||A|||$ is the modulus of the largest eigenvalue of $A$, while $||A||_{HS}$ is the quadratic mean of the eigenvalues; it is standard that there is $D_1>0$ such that 
$$\frac{1}{D_1}||A||_{HS}\le|||A|||\le D_1||A||_{HS}     \eqno (2.1)$$
for all $A$ symmetric. 

We define $M^d$ as the space of $d\times d$ symmetric matrices; we recall that 
$A\in M^d$ is positive semidefinite if 
$$(v,Av)\ge 0\qquad\forall v\in\R^d  .  \eqno (2.2)$$
It is standard that $B\in M^d$ is positive semidefinite if and only if 
$$(A,B)_{HS}\ge 0  \eqno (2.3)$$
for every $A\in M^d$ satisfying (2.2). We briefly prove this fact. Let $B$ satisfy (2.3); if we let let $A$ vary among the one-dimensional projections we easily see that all the eigenvalues of $B$ are positive, and (2.2) follows. For the converse, we note that, since $A$ and $B$ are symmetric and semi positive definite, their square roots are symmetric and real, which implies the first and third equalities below; for the second one, we use the fact that the trace of a product is invariant under cyclic permutations. 
$$tr(\tr AB)=tr(\sqrt{\tr A}\sqrt{\tr A}\sqrt B\sqrt B)=$$
$$tr(\sqrt B\sqrt{\tr A}\sqrt{\tr A}\sqrt B)=
tr(\tr(\sqrt{\tr A}\sqrt B)(\sqrt{\tr A}\sqrt B))\ge 0  .  $$
Essentially, (2.3) implies that the angle at the vertex of the cone of positive semidefinite matrices is smaller than $\frac{\pi}{2}$. 

An immediate consequence of (2.3) is the following: if $A,B,C\in M^d$, if $A\ge 0$ and $B\le C$, then 
$$(A,B)_{HS}\le(A,C)_{HS}  .  \eqno (2.4)  $$

Let now $(E,\hat d)$ be a compact metric space; in the following, $(E,\hat d)$ will be either one of $(G,\tilde d)$ or $(\Sigma,d_\g)$. We define $C(E,M^d)$ as the space of continuous functions from $E$ to $M^d$; for $A\in C(E, M^d)$ we define
$$||A||_\infty=\sup_{x\in E}||A(x)||_{HS}  .  $$
Let us call $\mc(E,M^d)$ the space of the Borel measures on $E$ valued in $M^d$. Putting on $M^d$ the Hilbert-Schmidt norm, we can define in the usual way the total variation $||\tau||$ of a measure $\tau\in\mc(E,M^d)$; clearly, $||\tau||$ is a scalar-valued, non-negative, finite measure on the Borel sets of $E$. 

If $\tau\in\mc(E,M^d)$ and $\fun{B}{E}{M^d}$ belongs to $L^1(E,||\tau||)$, we define the real number 
$$\int_E(B_x,\dr\tau(x))_{HS}\colon=
\int_E(B_x,T_x)_{HS}\dr||\tau||(x)$$
where $\tau=T_x||\tau||$ is the polar decomposition of $\tau$; we recall that 
$||T_x||_{HS}=1$ for $||\tau||$-a. e. $x\in S$. 

Several other products are possible; for instance, if $\fun{u,v}{E}{\R^d}$ are Borel vector fields such that 
$$||u(x)||\cdot||v(x)||\in L^1(E,||\tau||),$$
we can define the real number  
$$\int_E(u(x),\dr\tau(x)v(x))\colon=
\int_E(u(x),T_xv(x))\dr||\tau||(x)$$
where, again, $\tau=T_x||\tau||$ is the polar decomposition of $\tau$. Analogously, if 
$\fun{A}{E}{M^d}$ is a Borel field of matrices such that $||A_x||_{HS}\in L^1(E,||\tau||)$, we can define the two matrices 
$$\int_E A_x\dr\tau(x)\colon=
\int_E A_xT_x\dr||\tau||(x)$$
and 
$$\int_E\dr\tau(x)A_x\colon=
\int_ET_xA_x\dr||\tau||(x)  .  $$
If $Q\in C(E,M^d)$ and $\tau\in\mc(E,M^d)$, we define the scalar measure 
$(Q,\tau)_{HS}$ in the following way: if $B\subset E$ is Borel, then 
$$(Q,\tau)_{HS}(B)\colon=\int_B(Q,\dr\tau)_{HS} .  $$
In other words, $(Q,\tau)_{HS}=(Q_x,T_x)_{HS}||\tau||$.

By Riesz's representation theorem, $\mc(E,M^d)$ is the dual space of $C(E,M^d)$; the duality coupling 
$$\fun{\inn{\cdot}{\cdot}}{C(E,M^d)\times \mc(E,M^d)}{\R}   $$
is given by 
$$\inn{B}{\tau}=\int_E(B_x,\dr\tau(x))_{HS}  .  $$
By Lusin's theorem we get in the usual way that, if $B\subset E$ is a Borel set, then
$$||\tau||(B)=\sup\int_B(A,\dr\tau)_{HS}   \eqno (2.5)$$
where the $\sup$ is over all $A\in C(E,M^d)$ such that $||A||_\infty\le 1$. 

\vskip 1pc

\noindent{\bf Remark.} In the discussion above, we should have distinguished between 
$M^d$ and its dual $(M^d)^\ast$; strictly speaking, the dual of $C(E,M^d)$ is 
$\mc(E,(M^d)^\ast)$. In order to have a simpler notation, we identify $M^d$ and 
$(M^d)^\ast$ thanks to the Riemannian structure on $\R^d$. For the same reason, if 
$f\in C^1(\R^d,\R)$, we shall deal with its gradient $\nabla f$ and not with its differential 
$\dr f$.

\vskip 1pc

We shall say that $\tau\in\mc^+(E,M^d)$ if $\tau\in\mc(E,M^d)$ and $\tau(B)$ is a non-negative definite matrix for all Borel sets $B\subset E$. By Lusin's theorem, this is equivalent to 
$$\int_E(v_x,\dr\tau(x)v_x)\ge 0\qquad
\forall v\in C(E,\R^d)   .   $$
In turn, by (2.3) this is equivalent to 
$$\int_E(A_x,\dr\tau(x))_{HS}\ge 0  \eqno (2.6)$$
for all $A\in C(E,M^d)$ such that $A_x$ is positive semidefinite for all $x\in E$. 

Let now $Q\in C(E,M^d)$ such that $Q_x$ is positive-definite for all $x\in E$; since $E$ is compact there is $D_2>0$ such that 
$$\frac{1}{D_2}Id\le Q_x\le D_2Id\qquad
\forall x\in E  .  \eqno (2.7)$$
For $Q$ satisfying (2.7) we define $\pc_Q(E,M^d)$ as the set of all $\tau\in\mc^+(E,M^d)$ such that
$$\int_E(Q,\dr\tau)_{HS}=1  .  $$

\lem{2.1} Let $Q\in C(E,M^d)$ satisfy (2.7). Then, there is $D_3>0$ (depending on the constant $D_2$ of (2.7)) such that for all $\tau\in\mc^+(E,M^d)$ and all Borel sets 
$B\subset E$ we have 
$$||\tau||(B)\le D_3\cdot (Q,\tau)_{HS}(B)   \eqno (2.8)$$
where $||\cdot||$ denotes total variation. As a consequence, $\pc_Q(E,M^d)$ is a convex set of $\mc(E,M^d)$, compact for the weak$\ast$ topology. 

\proof By the definition of total variation, we must find $D_3>0$ with the following property: for all Borel sets $B\subset E$ and all countable Borel partitions 
$\{ B_i \}_{i\ge 1}$ of $B$ we have that 
$$\sum_{i\ge 1}||\tau(B_i)||_{HS}\le D_3\cdot (Q,\tau)_{HS}(B)  .  $$
By (2.1), this follows if we show that, for the constant $D_2$ of (2.7), 
$$\sum_{i\ge 1}|||\tau(B_i)|||\le D_2\cdot (Q,\tau)_{HS}(B) .  \eqno (2.9)$$
By the definition of $|||\tau(B_i)|||$ before (2.1), we can find unit vectors $v_i$ such that
$$|||\tau(B_i)|||=(v_i,\tau(B_i)v_i)\qquad\forall i\ge 1  .  $$
Let now $v\in\R^d$; the inequality below follows in a standard way from the fact that 
$\tau(B_i)$ is symmetric and non-negative-definite; the equality comes from the definition of the Hilbert-Schmidt product. 
$$(v,\tau(B_i)v)\le{\rm tr}(\tau(B_i))||v||^2=
(\tau(B_i),Id)_{HS}||v||^2  .  $$
Since $v_i$ has unit length, the last two formulas imply the first inequality below; the first equality follows since $\tau$ is a measure and $\{ B_i \}_{i\ge 1}$ is a partition of $B$. Since $\tau\in\mc^+(G,M^d)$, $\tau(B)$ is positive semidefinite; in particular, (2.4) holds and together with (2.7) implies the second inequality below. The last equality follows from the definition of the measure 
$(Q,\tau)_{HS}$.  
$$\sum_{i\ge 1}|||\tau(B_i)|||\le
\sum_{i\ge 1}(\tau(B_i),Id)_{HS}=
(\tau(B),Id)_{HS}\le
D_2\cdot\int_B (Q,\dr\tau)_{HS}=D_2\cdot (Q,\tau)_{HS}(B)   .  $$
This is (2.9) and we are done.

In order to prove the last assertion, we note that, by (2.6), $\mc^+(E,M^d)$ is a convex set of $\mc(E,M^d)$ clesed for the weak$\ast$ topology; as a consequence, also 
$\pc_Q(E,M^d)$ is a closed convex set, while (2.8) implies that it is relatively compact for the weak$\ast$ topology. 

\fin

\noindent{\bf Definitions.} Let $(E,\hat d)$ be a compact metric space with 
${\rm diam}(E)\le 1$. We define $C_+$ as the set of all the $A\in C(E,M^d)$ such that 
$A_x$ is positive-definite for all $x\in E$; since $E$ is compact, if $A\in C_+$, there is 
$\e>0$ (depending on $A$) such that
$$A_x\ge\e||A||_\infty Id\qquad\forall x\in E  .  \eqno (2.10)$$
For $a>0$ and $\nu\in(0,1]$ we define $C_+(E,a,\nu)$ as the set of all the $A\in C_+$ such that  
$$A_x e^{-a\hat d(x,y)^\nu}\le A_y\le A_x e^{a\hat d(x,y)^\nu}\qquad\forall x,y\in E .  \eqno (2.11)$$
We also define $C^\nu(E,M^d)$ as the set of all $\nu$-H\"older maps from $E$ to $M^d$, with the seminorm
$$||A||_\nu=\sup_{x\not=y\in E}
\frac{||A_x-A_y||_{HS}}{\hat d(x,y)^\nu}  .  $$
As lemma 2.3 below shows, the last two formulas are two different ways to look at the same seminorm, but we shall need both. 

\lem{2.2} Let $\e>0$ and let $A,B\in M^d$ such that 
$$A,B\ge\e Id  .  \eqno (2.12)$$
Then, there is $D_3=D_3(\e,B)>0$ such that 
$$Be^{-D_3||B-A||_{HS}}\le A\le Be^{D_3||B-A||_{HS}} .  \eqno (2.13)$$
For fixed $\e$, the function $D_3(\e,B)$ is bounded when $B$ is bounded. 

As a converse, there is $D_5>0$ such that the following holds. Let $A,B\in M^d$ be semi-positive definite and let us suppose that there is $D_4>0$ such that 
$$e^{-D_4}B\le A\le e^{D_4} B  .  \eqno (2.14)$$
Then, 
$$||B-A||_{HS}\le D_5(e^{D_4}-1)||A||_{HS}. \eqno (2.15)$$

\proof We begin with the direct part. Let $C\in M^d$; it is easy to see (for instance, choosing a base in which $C$ is diagonal) that   
$$C\le ||C||_{HS}Id  .  $$
Let $A,B\in M^d$; if we apply the formula above to $C=A-B$ we get that  
$$A\le B+||B-A||_{HS}Id  .  \eqno (2.16)$$
Since $B$ satisfies (2.12), this implies the first inequality below. 
$$A\le B\left(
1+\frac{1}{\e}||B-A||_{HS}
\right)  \le
Be^{
\frac{1}{\e}||B-A||_{HS}
}  .  $$
This yields the inequality on the right of (2.13). We prove the inequality on the left; the first inequality below is (2.16) with the names changed, the second one follows from (2.12).  
$$A\ge B-||B-A||_{HS}Id\ge 
B\left(
1-\frac{1}{\e}||B-A||_{HS}
\right)  .  $$
Again by (2.12) this implies that 
$$A\ge\left\{
\eqalign{
B\left( 1-\frac{1}{\e}||B-A||_{HS} \right) &\txt{if}||B-A||_{HS}\le\frac{\e}{2}\cr
\e Id&\txt{if}||B-A||>\frac{\e}{2}  .
}
\right.     $$
The left hand side of (2.13) now follows from two facts: the first one is that, for $D_3$ large enough, 
$$1-\frac{t}{\e}\ge e^{-D_3 t}\txt{if}0\le t\le\frac{\e}{2}  .  $$
The second one is the formula below. The first inequality comes taking $||B-A||_{HS}\ge\e$, the second one taking $\g>0$ so small that $\g B\le Id$; the third one taking $D_3$ so large that 
$$\frac{1}{\g}e^{-\frac{\e}{2}D_3}\le\e  .  $$
$$B e^{-D_3||B-A||_{HS}}\le
B e^{-D_3\frac{\e}{2}}\le
\frac{1}{\g}Id\cdot e^{-D_3\frac{\e}{2}}\le
\e Id  .  $$

We prove the converse. By (2.14) we have that 
$$-(1-e^{-D_4})(Bx,x)\le ((A-B)x,x)\le (e^{D_4}-1)(Bx,x)  .  $$
By the definition of $|||\cdot|||$ and the fact that $e^{D_4}-1\ge 1-e^{-D_4}$ this implies that 
$$ |||A-B|||\le(e^{D_4}-1)|||B|||   .  $$
Now (2.15) follows from (2.1). 

\fin

\lem{2.3} Let $a>0$ and let $\nu\in(0,1]$; then, the following holds. 

\noindent 1) The sets $C_+$ and $C_+(E,a,\nu)$ are convex cones in $C(E,M^d)$ which satisfy (1.1). 

\noindent 2) There is $D_6>0$ such that for all 
$A\in C_{+}(E,a,\nu)$ we have that 
$$||A_x-A_y||\le D_6||A||_\infty\cdot\hat d(x,y)^\nu\qquad\forall x,y\in G  .  \eqno (2.17)$$ 

Conversely, if  $A\in C_{+}\cap C^\nu(E,M^d)$, then $A\in C_+(E,a,\nu)$ for some $a>0$ (which depends on $A$).

\proof We don't dwell on the proof of point 1), since it follows immediately from the definitions of $C_+$ and $C_+(a,\nu)$. 

We prove point 2). Let $A\in C_+(E,a,\nu)$; this means that, if $x,y\in E$, 
$$e^{-a\hat d(x,y)^\nu}A_y\le A_x\le e^{a\hat d(x,y)^\nu}A_y   .   $$
This is (2.14) for $D_4=a\hat d(x,y)^\nu$; by lemma 2.2, (2.15) holds, i. e. 
$$||A_x-A_y||_{HS}\le D_5   \left(
e^{a\hat d(x,y)^\nu} -1
\right)    ||A_x||_{HS} \qquad\forall x,y\in E  .  $$
Since we are supposing that the diameter of $E$ is 1 (for $E=G$ this is (1.9)), we get (2.17).  

Conversely, let $A\in C_+\cap C^\nu(E,M^d)$; since $E$ is compact, we easily see that there is $\e>0$ such that $A_x\ge\e Id$ for all $x\in E$. Thus, setting $A=A_x$ and 
$B=A_y$, we have that $A$ and $B$ satisfy (2.12); by lemma 2.2 also (2.13) holds, i. e. 
$$A_y e^{-D_3||A_x-A_y||_{HS}}\le A_x\le 
A_y e^{D_3||A_x-A_y||_{HS}}     $$
for some $D_3=D_3(\e,A_y)>0$. Since $A\in C^\nu(E,M^d)$, we have that 
$$||A_x-A_y||_{HS}\le D_7\hat d(x-y)^\nu    $$
The last two formulas imply that 
$$A_ye^{
-D_7 D_3(\e,A_y)\hat d(x,y)^\nu
}  \le A_x\le
A_ye^{
D_7 D_3(\e,A_y)\hat d(x,y)^\nu
}  .   $$
Thus, $A\in C_+(E,a,\nu)$ if 
$$a\ge D_7\sup\{ D_3(\e,A_y)\st y\in E  \} . $$
Note that the term on the right is finite since, by lemma 2.2, $D_3(\e,\cdot)$ is bounded on bounded sets and $||A||_\infty$ is finite.

\fin

\vskip 2pc
\centerline{\bf \S 3}
\centerline{\bf The Ruelle operator}
\vskip 1pc

From now on, we suppose that (F1)-(F4) hold; we let $(G,\tilde d)$ be the fractal defined by (1.8) with the distance $\tilde d$ induced by the immersion in $\R^d$. 

We define the Ruelle operator $\L_G$ on $G$ as 
$$\fun{\L_G}{C(G,M^d)}{C(G,M^d)}  $$
$$(\L_G A)(x)=
\sum_{i=1}^n
\tr D\psi_i(x)  A_{\psi_i(x)}  D\psi_i(x)   .  \eqno (3.1)$$

We also define a Ruelle operator $\L_\Sigma$ on $C(\Sigma,M^k)$: first, if $x=(x_0x_1\dots)\in\Sigma$, we set $(ix)=(ix_0x_1\dots)$. Then, we define  
$$\fun{\L_\Sigma}{C(\Sigma,M^k)}{C(\Sigma,M^k)}$$
$$(\L_\Sigma A)(x)=\sum_{i=1}^n 
\tr D\psi_i|_{\Phi(x)}A_{(ix)} D\psi_i|_{\Phi(x)} .  \eqno (3.2)$$
Note also that $\fun{}{x}{D\psi_i|_{\Phi(x)}}$ is $\nu_0$-H\"older, since $D\psi_i$ is 
$\nu_0$-H\"older by (F1) and we saw in section 1 that $\Phi$ is Lipschitz; if $A=\tilde A\circ\Phi$ for some $\tilde A\in C(G,M^d)$,  we get by the second formula of (1.15) that 
$$\L_\Sigma A(x)= 
\sum_{i=1}^n 
\tr D\psi_i|_{\Phi(x)}\tilde A_{\psi_i\circ\Phi(x)} D\psi_i|_{\Phi(x)}  .  $$

The next lemma shows that the fixed points of the adjoint of $\L_G$, which we call 
$\L^\ast_G$, induce a self-similar form on $C^1(\R^d)$.

\lem{3.1}  Let $G$ be  a fractal satisfying (F1)-(F4).  Let 
$(\b,\tau)\in(0,+\infty)\times\mc^+(G,M^d)$ be such that 
$$\L_G^\ast\tau=\b\tau   .    $$
Let $f,g\in C^1(\R^d,\R)$ and let us define, with the notation of section 2 for the integral, 
$$\ec_\tau(f,g)=
\int_G(\nabla f(x),\dr\tau\nabla g(x)) . $$
Then, 
$$\sum_{i=1}^n\ec_\tau(f\circ\psi_i,g\circ\psi_i)=
\b\ec_\tau(f,g)   .  \eqno (3.3)$$

\proof By the polarisation identity, it suffices to show (3.3) when $f=g$; we shall take advantage of the fact that $\nabla f\otimes\nabla f\in C(G,M^d)$, i. e. is in the domain of 
$\L_G$. We recall that, if $a\in\R^d$ and $A\in M^d$, then 
$$(a,Aa)=(a\otimes a,A)_{HS}     $$
where by $a\otimes a$ we denote the tensor product of the column vector $a$ with itself:
$$a\otimes a=a\cdot\tr a  .  $$
As we shall see in the formula below, this explains the position of the transpose sign in (3.1). 

The definition of $\ec_\tau$ and the formula above imply the first equality below; the second one comes from the chain rule (recall that 
$\nabla(f\circ\psi_i)=\tr D\psi_i\cdot\nabla f$) and the definition of $\L_G$; this third one follows since 
$\L^\ast_G\tau=\b\tau$ and the last one is again the definition of $\ec_\tau$. 
$$\sum_{i=1}^n\ec_\tau(f\circ\psi_i,f\circ\psi_i)=
\sum_{i=1}^n\int_G(\nabla(f\circ\psi_i)(x)\otimes\nabla(f\circ\psi_i)(x),\dr\tau(x))_{HS}=$$
$$\int_G(\L_G(\nabla f\otimes\nabla f)(x),\dr\tau(x))_{HS}=$$
$$\b\int_G(\nabla f\otimes\nabla f,\dr\tau)_{HS}=
\b\ec_\tau(f,g)   .  $$

\fin 

\noindent{\bf Remark.} We can read lemma 3.1 as a statement about the push-forward of the measure $\tau$, the positive eigenvector of $\L^\ast$. Indeed, let $f,g\in C^1(\R^d)$; by (F3), $f\circ F$ and $g\circ F$ are not defined at the points of $\fc$, which are a finite set. As we shall see in lemma 5.3 below, the measure $\tau_G$ of point 2) of theorem 1 is non-atomic; in particular, the points where $f\circ F$ and $g\circ F$ are not defined have measure zero. Together with (1.13), this implies the first equality below; the second one follows from lemma 3.1; the last one follows from the chain rule.
$$\int_G(\nabla f,\dr\tau\nabla g)=
\frac{1}{n}\sum_{i=1}^n
\int_G(\nabla (f\circ F\circ\psi_i),\dr\tau\nabla (g\circ F\circ\psi_i))=$$
$$\frac{\b}{n}\int_G(\nabla(f\circ F),\dr\tau\nabla(g\circ F))=
\frac{\b}{n}\int_G(
\tr DF(x)\nabla f|_{F(x)}\otimes\tr DF(x)\nabla g|_{F(x)},\dr\tau(x)
)_{HS}   .    \eqno (3.4)$$
In other words, denoting by $F_\ast\omega$ the pull-back by $F$ of the two-tensor 
$\omega$, we have that 
$$\tau= \frac{\b}{n}(F_\ast)_\sharp \tau $$
where the "push-forward" $(F_\sharp)_\sharp\tau$ is defined as in the last term on the right in the formula above. Though this is not the standard push-forward operator, it is natural if we regard $\tau$ not as a measure, but as a linear operator on 2-tensors.

\vskip 2pc
\centerline{\bf \S 4}
\centerline{\bf Fixed points of the Ruelle operator}
\vskip 1pc

We shall suppose that the maps $\{ \psi_i \}_{i=1}^n$ satisfy the following nondegeneracy condition; it is stronger than the one in [12], but it allows us to use the Perron-Frobenius theorem without modifications. 

\vskip 1pc

\noindent{\bf (ND)} We suppose that, for all $v\in\R^d\setminus\{ 0 \}$ and all $x\in G$, the set $\{ \tr D\psi_i(x)v \}_{i=1}^n$ generates $\R^d$. Actually, we ask for a quantitative version of this, i. e. that there is $b>0$ such that the following holds. Let 
$v,v_0\in\R^d$ and let $x\in G$; then, there is $\bar i\in(1,\dots,n)$, depending on $x$, $v$ and $v_0$,  such that 
$$( D\psi_{\bar i}(x)v,v_0)\ge b||v||\cdot||v_0||  .  \eqno (4.1)$$
If we denote by $P$ the orthogonal projection on $v_0$, the formula above implies the inequality below. 
$$||P\cdot D\psi_{\bar i}(x)v||\ge  b||v||  .  \eqno (4.2)$$

\vskip 1pc

It is easy to verify that the harmonic Sierpinski gasket of section 1 satisfies (ND); as we shall see in the next lemma, (ND) implies a bound from below on $\L A$.

\lem{4.1} Let the maps $\{ \psi_i \}_{i=1}^n$ satisfy (F1)-(F4) and (ND). Let $(E, \L)$ denote either one of $(G,\L_G)$ or $(\Sigma,\L_\Sigma)$. 

Then, for all $a>0$ there is $D_1=D_1(a,b)>0$ such that the following happens. Let 
$\nu\in(0,\nu_0]$ and let $A\in{\cal C}_+(E,a,\nu)$; then, for all $x\in E$, 
$$\frac{1}{D_1(a,b)}||A||_\infty\cdot Id\le
(\L A)_x\le
D_1(a,b)||A||_\infty\cdot Id  .  \eqno (4.3)$$

\proof We prove the left hand side of (4.3); for the right hand side it suffices to note that 
$\L$ is continuous from the $||\cdot||_\infty$ topology to itself. 

By compactness, there is $x_{max}\in E$ such that 
$$||A_{x_{max}}||_{HS}=||A||_{\infty}   .   \eqno (4.4)$$
By the definition of $|||A_{x_{max}}|||$ we can find $v_{max}\in\R^d$ with $||v_{max}||=1$ such that
$$|||A_{x_{max}}|||=(A_{x_{max}}v_{max},v_{max})   .  \eqno (4.5)$$
By (2.11) and the fact that ${\rm diam}(E)=1$ we have that 
$$e^{-a}(A_{x_{max}}v,v)\le
(A_xv,v)\le
e^a (A_{x_{max}}v,v)\qquad
\forall x\in E,\quad\forall v\in\R^d  .  $$
Let $v\in\R^d$ and let $\bar i\in(1,\dots,n)$; the formula above implies the first inequality below. Since $v_{max}$ is an eigenvector of the symmetric matrix $A_{x_{max}}$, we have that $A_{x_{max}}$ preserves the space generated by $v_{max}$ and its orthogonal complement. Thus, if we denote by $P$ the orthogonal projection on $v_{max}$, we get the second inequality below. Next, we choose $\bar i$ in such a way that (4.1) holds with  $v_0=v_{max}$; the choice of $\bar i$ depends on $x$ and $v$. By (4.2) we get the third inequality below; the last equality comes from (4.5) and the last inequality from (2.1) and (4.4). 
$$(A_{\psi_{\bar i}(x)} D\psi_{\bar i}(x)v, D\psi_{\bar i}(x)v)\ge$$
$$e^{-a}(A_{x_{max}} D\psi_{\bar i}(x)v, D\psi_{\bar i}(x)v)\ge
e^{-a}(A_{x_{max}}P D\psi_{\bar i}(x)v,P D\psi_{\bar i}(x)v)\ge$$
$$e^{-a} b^2(A_{x_{max}}v_{max},v_{max})\cdot||v||^2=
e^{-a}b^2 |||A_{x_{max}}|||\cdot ||v||^2\ge
D_4\cdot e^{-a} b^2\cdot ||A||_{\infty}\cdot||v||^2  .  $$
We choose $\bar i$ as above; the definition of $\L$ implies the first inequality below; the second one comes from the formula above.  
$$((\L A)_xv,v)\ge 
(A_{\psi_{\bar i}(x)} D\psi_{\bar i}(x)v,  D\psi_{\bar i}(x)v)\ge$$
$$D_4\cdot e^{-a} b^2\cdot ||A||_{\infty}\cdot ||v||^2
\qquad\forall x\in E  .  $$

\fin

\lem{4.2} Let the maps $\{ \psi_i \}_{i=1}^n$ satisfy (F1)-(F4) and (ND) for some $b>0$. Then,  for all $a>0$ there is $\o(a,b)>0$ such that, if $A\in C_+(G,a,\nu_0)$ and 
$$||D\psi_i||_{\nu_0}\le\o(a,b)\qquad \forall i\in(1,\dots,n)  ,  \eqno (4.6)$$
then for all 
$x,y\in G$ we have that 
$$e^{
-||x-y||^{\nu_0}
}
\sum_{i=1}^n\tr D\psi_i(x)A_{\psi_i(x)}  D\psi_i(x)\le$$
$$\sum_{i=1}^n\tr D\psi_i(y)A_{\psi_i(x)}  D\psi_i(y)\le$$
$$e^{
||x-y||^{\nu_0}
}
\sum_{i=1}^n\tr D\psi_i(x)A_{\psi_i(x)}  D\psi_i(x)    .    \eqno (4.7)_G$$
Note that, if the maps $\psi_i$ are affine and satisfy (ND), (4.6) is always verified. 

Analogously, possibly reducing $\o(a,b)$ in (4.6), for all $A\in C_+(\Sigma,a,\nu_0)$ and all $x,y\in \Sigma$ we have that 
$$e^{
-d_\g(x,y)^{\nu_0}
}
\sum_{i=1}^n\tr D\psi_i|_{\Phi(x)}A_{(ix)}  D\psi_i|_{\Phi(x)}\le$$
$$\sum_{i=1}^n\tr D\psi_i|_{\Phi(y)}A_{(ix)}  D\psi_i|_{\Phi(y)}\le$$
$$e^{
d_\g(x,y)^{\nu_0}
}
\sum_{i=1}^n\tr D\psi_i|_{\Phi(x)}A_{(ix)}  D\psi_i|_{\Phi(x)}    .    \eqno (4.7)_\Sigma$$

\proof We shall prove $(4.7)_G$, since $(4.7)_\Sigma$ is analogous. We begin recalling an inequality on matrices. Let $B\in M^d$ be positive semidefinite and let $C,C^\prime$ two invertible matrices; we suppose that 
$$||C||_{HS},  ||C^\prime||_{HS}\le D_1      \eqno (4.8)$$
for some $D_1>0$.

It is easy to see that, since $B$ is symmetric, 
$$(\tr C^\prime BCx,x)=(\tr CBC^\prime x,x)  .  $$
Together with a simple calculation, this implies that 
$$(\tr C^\prime BC^\prime x,x)=(BC^\prime x,C^\prime x)=$$
$$(BCx,Cx)+2(B(C^\prime-C)x,Cx)+
(B(C^\prime-C)x,(C^\prime-C)x)  .   $$
Since $||C-C^\prime||_{HS}\le 2D_1$ by (4.8), this implies that, for some $D_2>0$ depending only on $D_1$, but not on $C,C^\prime$ and $B$, 
$$\tr CBC-D_2||B||_{HS}\cdot||C-C^\prime||_{HS}\cdot Id\le
\tr C^\prime B C^\prime\le
\tr CBC+D_2||B||_{HS}\cdot||C-C^\prime||_{HS}\cdot Id   .   $$
We set 
$$B=A_{\psi_i(x)},\quad    C^\prime= D\psi_i(y)
\txt{and} C= D\psi_i(x)    $$
which immediately implies that
$$||C-C^\prime||_{HS}\le ||D\psi_i||_{\nu_0}\cdot ||x-y||^{\nu_0}  .  $$
Recalling that by (F1) (4.8) holds with $D_1=\eta$, we get from the last three formulas that 
$$\tr D\psi_i(x) A_{\psi_i(x)} D\psi_i(x)-
D_5|| A_{\psi_i(x)}  ||_{HS}\cdot ||D\psi_i||_{\nu_0}\cdot
||x-y||^{\nu_0}Id\le$$
$$\tr D\psi_i(y) A_{\psi_i(x)} D\psi_i(y)\le$$
$$\tr D\psi_i(x) A_{\psi_i(x)} D\psi_i(x)+
D_5|| A_{\psi_i(x)} ||_{HS}\cdot ||D\psi_i||_{\nu_0}\cdot
||x-y||^{\nu_0}Id    .    $$
Summing over $i\in(1,\dots,n)$ and setting
$$||D\psi||_{\nu_0}\colon=
\sup_{i\in (1,\dots,n)}||D\psi_i||_{\nu_0}$$
we get
$$\sum_{i=1}^n\tr D\psi_i(x) A_{\psi_i(x)} D\psi_i(x)-
D_5||D\psi||_{\nu_0}\cdot\sum_{i=1}^n|| A_{\psi_i(x)} ||_{HS}\cdot
||x-y||^{\nu_0}Id\le$$
$$\sum_{i=1}^n\tr D\psi_i(y) A_{\psi_i(x)} D\psi_i(y)\le$$
$$\sum_{i=1}^n\tr D\psi_i(x) A_{\psi_i(x)} D\psi(x)+
||D\psi||_{\nu_0}\cdot D_5\sum_{i=1}^n|| A_{\psi_i(x)} ||_{HS}\cdot
||x-y||^{\nu_0}Id    .    $$
Since $A\in C_+(a,\nu)$ we can apply lemma 4.1 and get that there is 
$D_6=D_6(a,b)$ such that 
$$[1-D_6(a,b)||D\psi||_{\nu_0}\cdot||x-y||^{\nu_0}]\sum_{i=1}^n\tr D\psi_i(x)A_{\psi_i(x)}
D\psi_i(x)\le$$
$$\sum_{i=1}^n\tr D\psi_i(y)A_{\psi_i(x)} D\psi_i(y)\le$$
$$[1+D_6(a,b)||D\psi||_{\nu_0}\cdot||x-y||^{\nu_0}]
\sum_{i=1}^n\tr D\psi_i(x)A_{\psi_i(x)} D\psi_i(x)  .  $$
We take 
$$\o(a,b)=\frac{1}{4D_6(a,b)}    $$
and we recall that, if $x\in[0,1]$, 
$$1-\frac{1}{4} x\ge e^{-x}
\txt{and}
1+\frac{1}{4} x\le e^{x}  .  $$
From the last three formulas, (4.6) and the fact that ${\rm diam}(E)=1$ we get that  
$$e^{
-||x-y||^{\nu_0}
}
\sum_{i=1}^n\tr D\psi_i(x)A_{\psi_i(x)}  D\psi_i(x)\le$$
$$\sum_{i=1}^n\tr D\psi_i(y)A_{\psi_i(x)}  D\psi_i(y)\le$$
$$e^{
||x-y||^{\nu_0}
}
\sum_{i=1}^n\tr D\psi_i(x)A_{\psi_i(x)}  D\psi_i(x)        $$
which is $(4.7)_G$.

\fin

\lem{4.3} Let $(E,\L)$ be either one of $(G,\L_G)$ or $(\Sigma,\L_\Sigma)$; let (F1)-(F4) and (ND) hold. Then, there is $a_0>0$ such that, for $a>a_0$ and 
$||D\psi_i||_{\nu_0}\le\o(a,b)$, 
$$\L(C_+(E,a,\nu_0))\subset C_+(E,a-1,\nu_0)   .  \eqno (4.9)$$

\proof We follow [19].  It is immediate from the definition of $\L$ that 
$\L(C_+)\subset C_{+}$. Thus, it suffices to show that, if $A\in C_+$ satisfies (2.11) for $a$ and $\nu_0$, then $\L A$ satisfies (2.11) for $a-1$ and $\nu_0$, provided $a$ is large enough. 

We shall prove the lemma on $G$, since the proof on $\Sigma$ is analogous. Let 
$x,y\in G$; the first equality below is the definition of $\L_G$ in (3.1); the first inequality is the left hand side of $(4.7)_G$ and holds if $||D\psi_i||_{\nu_0}\le\o(a,b)$. The second one follows from two facts: the map $\fun{}{A}{\tr BAB}$ is order-preserving and 
$A\in C_{+}(a,\nu_0)$, i. e. it satisfies (2.11). The third inequality comes from the fact that $\eta$ is the common Lipschitz constant of the maps $\psi_i$, i. e. formula (1.7). 
$$(\L_G A)(y)=
\sum_{i=1}^n\tr D\psi_i(y) A_{\psi_i(y)}  D\psi_i(y)\le$$
$$e^{ ||x-y||^{\nu_0}}
\sum_{i=1}^n
\tr D\psi_i(x) A_{\psi_i(y)}  D\psi_i(x)  \le$$
$$e^{
 ||x-y||^{\nu_0}+a ||\psi_i(x)-\psi_i(y)||^{\nu_0}
}  
\sum_{i=1}^n\tr D\psi_i(x) A_{\psi_i(x)}  D\psi_i(x)\le$$
$$e^{(1 +a\eta^{\nu_0})||x-y||^{\nu_0}}
\sum_{i=1}^n\tr D\psi_i(x) A_{\psi_i(x)}  D\psi_i(x)  .   $$
Since $\eta\in (0,1)$ we can choose $a_0$ so large that, for $a\ge a_0$, 
$$1+a\eta^{\nu_0}\le a-1  .$$
From the last two formulas we get that 
$$(\L_G A)(y)\le e^{(a-1)||x-y||^{\nu_0}}(\L_G A)(x)  .  $$
The opposite inequality follows similarly, implying (4.9).

\fin 

\noindent{\bf Definitions.} Let $\l_1,\e>0$; we denote by $C^{\e}_+(E,\l_1a,\nu_0)$ the subset of the $A\in C_+(E,\l_1 a,\nu_0)$ such that, for all $x\in E$, 
$$A_x\ge\e ||A||_\infty Id   .  \eqno (4.10)$$

Moreover, we shall call $\theta^{(a,\nu_0)}$ the hyperbolic distance on $C_+(E,a,\nu_0)$ and 
$\theta_+$ the hyperbolic distance on $C_+$; we recall that the hyperbolic distance on a cone has been defined at the beginning of section 1. 

\vskip 1pc

\lem{4.4} Let $E$ be either one of $G$ or $\Sigma$; let $a>0$. Then, the following holds.

\noindent 1) $(\frac{C_+(E,a,\nu)}{\simeq},\theta^{(a,\nu)})$ is a complete metric space. 

\noindent 2) Let $C^{\e}_+(E,\l_1a,\nu_0)$ be defined as in (4.10). If $\l_1\in(0,1)$, then
$${\rm diam}_{\theta^{(a,\nu_0)}}C^{\e}_+(E,\l_1a,\nu_0)<+\infty  .  $$

\proof Again we follow closely [19]; we begin with point 1). We consider  a Cauchy sequence 
$\{ A_n \}_{n\ge 1}$ in $(\frac{C_+(E,a,\nu_0)}{\simeq},\theta^{(a,\nu_0)})$; we choose the representatives which satisfy
$$||A_n||_\infty=1\qquad\forall n\ge 1  .  \eqno (4.11)$$

\noindent{\bf Step 1.} We begin to show that $\{ A_n \}_{n\ge 1}$ converges uniformly to 
$A\in C(E,M^d)$. 

Since $C_+(E,a,\nu_0)\subset C_+$, the definition of the hyperbolic distance implies that 
$\theta_+\le\theta^{(a,\nu_0)}$; thus, $\{ A_n \}_{n\ge 1}$ is Cauchy also for the $\theta^+$ distance; the definition of $\theta^+$ in (1.4) implies that
$$\frac{\b^+(A_m,A_n)}{\a^+(A_m,A_n)}\tends 1
\txt{as} n,m\tends+\infty  .  \eqno (4.12)$$
By the definition of $\a^+$ and $\b^+$ in (1.2) and (1.3) respectively, we get that, for all 
$x\in E$, 
$$\a^+(A_m,A_n)A_m(x)\le A_n(x)\le\b^+(A_m,A_n)A_m(x)  .  $$
Let $\d>0$; by (4.12) we have that, for $n$ and $m$ large enough and all $x\in E$, 
$$(1-\d)\b^+(A_m,A_n)A_m(x)\le A_n(x)\le\b^+(A_m,A_n)A_m(x)  .  $$
By the converse part of lemma 2.2 this implies that, for $n$ and $m$ large, 
$$||A_n-\b^+(A_m,A_n)A_m||_\infty\le D_5
\cdot\d ||A_n||_\infty  .  $$
By (4.11) and the triangle inequality this implies that $\b^+(A_m,A_n)\tends 1$; again by the formula above, $\b^+(A_m,A_n)\tends 1$ implies that $\{ A_n \}_{n\ge 1}$ is Cauchy for $||\cdot||_\infty$; thus, there is $A\in C(G,M^d)$ such that $A_n\tends A$ uniformly. 

\noindent{\bf Step 2.} We show that $\theta^{(a,\nu_0)}(A_n,A)\tends 0$ as 
$n\tends+\infty$. It is easy to see that $\theta^{(a,\nu_0)}$ is lower semicontinuous under uniform convergence; together with step 1, this yields the first inequality below. Since 
$\{ A_n \}_{n\ge 1}$ is Cauchy for $\theta^{(a,\nu_0)}$, there is $\d_n\tends 0$ such that also the second inequality below holds.
$$\theta^{(a,\nu_0)}(A,A_n)\le
\liminf_{m\tends+\infty}\theta^{(a,\nu_0)}(A_m,A_n)\le\d_n  .  $$

\noindent{\bf End of the proof of point 1.} It only remains to prove that 
$A\in C_+(E,a,\nu_0)$. First of all, $A$ satisfies (2.11), since this condition is closed under uniform convergence. We have to show that $A_x$ is positive-definite for all $x\in E$. We recall that $A_{n,x}$ is positive-definite for all $x\in E$, since $A_n\in C_+(E,a,\nu_0)$; since $\theta^+(A_n,A)<+\infty$, we have that $\a^+(A_n,A)>0$; since by (1.2) 
$A\ge \a^+(A_n,A)A_n$ and $A_n$ satisfies (2.10) for some $\e>0$ (which depends on $A_n$), we get that $A_x$ is positive-definite for all $x\in E$. 

\noindent{\bf Proof of point 2).} The proof is in two steps: first, we show that 
$${\rm diam}_{\theta^+}C^{\e}_+(E,\l_1a,\nu_0)<+\infty  \eqno (4.13)$$
and then that
$${\rm diam}_{\theta^{(a,\nu_0)}}C^{\e}_+(E,\l_1a,\nu)\le
{\rm diam}_{\theta^+}C^{\e}_+(E,\l_1a,\nu_0)+D_5(\l_1)  .  \eqno (4.14)$$

\noindent{\bf Step 3.} We prove (4.13); this follows if we show that 
$C^{\e}_+(E,\l_1a,\nu_0)$ is compact in $(C_+,\theta^+)$. Thus, let 
$\{ A_n \}_{n\ge 1}\subset C^{\e}_+(E,\l_1a,\nu_0)$; we can suppose that 
$\{ A_n \}_{n\ge 1}$ is normalised, i. e. that (4.11) holds. Now, point 2) of lemma 2.3 implies that the H\"older seminorm of $A_n$ is bounded. Thus, by 
Ascoli-Arzel\`a\ there is a subsequence $\{ A_{n_h} \}_{h\ge 1}$ which converges uniformly to $A\in C(E,M^d)$; we see as in point 1) that $A$ satisfies (2.11); it also satisfies (4.10) because this formula is stable under uniform convergence. Thus, 
$A\in C^{\e}_+(E,\l_1a,\nu_0)$. Since $A_{n_h}$ and $A$ satisfy (4.10) and 
$A_{n_h}\tends A$ uniformly, we can apply the direct part of lemma 2.2 and get that 
$\theta^+(A_{n_h},A)\tends 0$, ending the proof of compactness. 

\noindent{\bf Step 4.} We prove (4.14) on the general space $E$ with distance $\hat d$. For starters, let us see how $\a^+$ and $\a^{(a,\nu_0)}$ are related. Let 
$A_1,A_2\in C_+(E,a,\nu_0)$; by the definition of $\a^{(a,\nu_0)}(A_1,A_2)$ we have that 
$$A_2-\a^{(a,\nu_0)}(A_1,A_2)A_1\in C_+(E,a,\nu_0)$$
which by (2.11) implies that, for all $x,y\in E$,  
$$e^{-a\hat d(x,y)^{\nu_0}}[
A_2(x)-\a^{(a,\nu_0)}(A_1,A_2)A_1(x)
]  \le
A_2(y)-\a^{(a,\nu_0)}(A_1,A_2)A_1(y)\le$$
$$e^{a\hat d(x,y)^\nu_0}[
A_2(x)-\a^{(a,\nu_0)}(A_1,A_2)A_1(x)
]  .  $$
Rearranging the terms of the inequality on the left, we get that 
$$e^{-a\hat d(x,y)^{\nu_0}}A_2(x)-A_2(y)\le
[
e^{-a\hat d(x,y)^{\nu_0}}A_1(x)-A_1(y)
]\a^{(a,\nu_0)}(A_1,A_2)  $$
for all $x,y\in E$. Since $A_1,A_2$ satisfy (2.11) for $\l_1a$, the formula above implies that
$$A_2(x)\le
\frac{
e^{-a\hat d(x,y)^{\nu_0}}-e^{-\l_1a d(x,y)^{\nu_0}}
}{
e^{-a\hat d(x,y)^{\nu_0}}-e^{\l_1a d(x,y)^{\nu_0}}
}    \a^{(a,\nu)}(A_1,A_2)A_1(x)      .   \eqno (4.15)$$
We set 
$$D_6=\sup\left\{
\frac{z-z^{\l_1}}{z-z^{-\l_1}}
\st z\in(0,1)
\right\}   $$
and by a function study we see that $D_6\in(0,1)$. By (4.15) we get that 
$$A_2(x)\le D_6\a^{(a,\nu_0)}(A_1,A_2)A_1(x)$$
which by the definition of $\a^+$ implies that
$$D_6\a^{(\a,\nu)}(A_1,A_2)\ge \a^+(A_1,A_2)  .  \eqno (4.16)$$
Analogously, we can set 
$$D_7=\inf\left\{
\frac{z-z^{-\l_1}}{z-z^{\l_1}}
\st z>1
\right\}   .  $$
A function study shows that $D_7>1$ and the same argument that yielded (4.16) yields 
$$\b^{(a,\nu)}(A_1,A_2)\le D_7\b^+(A_1,A_2)  .  $$
Using this, (4.16) and the definition of $\theta^+$ in (1.4) we get that
$$\theta^{(a,\nu)}(A_1,A_2)\le\theta^+(A_1,A_2)+
\log D_7-\log D_6  $$
which ends the proof of (2.18).

\fin

\prop{4.5} Let $(E,\L)$ be either one of $(G,\L_G)$ or $(\Sigma,\L_\Sigma)$. Let (F1)-(F4) and (ND) with constant $b>0$ hold. Then, if $\sup_{i\in(1,\dots,n)}||D\psi_i||_{\nu_0}$ is small enough, the following holds.

\noindent 1) There is a simple, positive eigenvalue $\b$ of
$$\fun{\L}{C(E,M^d)}{C(E,M^d)}  .  $$
Denoting by $Q$ the eigenfunction of $\b$, we have that $Q\in C_+(a,\nu_0)$. In particular, $Q_x$ is positive-definite for all $x\in E$. If the maps $\psi_i$ are affine, then 
there is $\bar Q\in M^d$ such that $Q_x=\bar Q$ for all $x\in E$. 

\noindent 2) Recall that after formula (2.7) we defined $\pc_Q(E,M^d)$; we assert that there is $\tau\in \pc_Q(E,M^d)$ such that $\L^\ast\tau=\b\tau$.

\noindent 3) If $B\in C(E, M^d)$, we have  
$$\frac{1}{\b^l}\L^l B\tends Q\int_E (B,\dr\tau)_{HS}  \eqno (4.17)$$
uniformly on $E$. Note that this implies that the measure $\tau$ of the previous point is unique and that $\b$ is simple eigenvalue of $\L$. Moreover, if $B\in C_+(E,a,\nu)$ with 
$\nu\in(0,\nu_0]$, the convergence above is exponentially fast. 

\proof {\bf Step 1.} Let $\sup_{i\in(1,\dots,n)}||D\psi_i||_{\nu_0}$ be so small that lemma 4.3 hold. Since $\fun{\L}{C(E,M)}{C(E,M)}$ is continuous, the left hand side of (4.3) implies that, possibly increasing the constant $D_1(a,b)$, 
$$(\L A)_x\ge\frac{1}{D_1(a,b)}||\L A||_\infty\cdot Id\qquad\forall x\in E  .  $$
Together with lemma 4.3 this implies that 
$$\L(C_+(E,a,\nu_0))\subset C^{\e}_+(E,\l_1a,\nu_0)  $$
for some $\l_1\in(0,1)$ and $\e$ equal to the constant $\frac{1}{D_1(a,b)}$ of the formula above; by point 2) of lemma 4.4 this implies that
$${\rm diam}_{\theta^{(a,\nu)}}\L(C_+(E,a,\nu_0))<+\infty  .  $$
By point 1) of proposition 1.1, we get that $\L$ is a contraction of $C_+(E,a,\nu)$ into itself; since 
$$\left(
\frac{C_+(E,a,\nu_0))}{\simeq},\theta^{(a,\nu_0)}
\right)    $$
is complete by point 1) of lemma 4.4, we get that $\L$ has a unique fixed point in 
$\frac{C_+(E,a,\nu_0))}{\simeq}$. In other words, there are 

\noindent 1) $Q\in C_+(E,a,\nu_0)$, unique up to multiplication by a scalar, and 

\noindent 2) a unique $\b\in(0,+\infty)$ such that 
$$\L Q=\b Q  .  \eqno (4.18)$$
Since $Q$ is unique up to multiplication by a scalar, we can normalise it in such a way that $||Q||_\infty=1$. 

\noindent{\bf Step 2.} If the maps $D\psi_i$ are constant, we see that, if $A$ is a constant matrix, then $\L A$ is constant too. Applying the Perron-Frobenius theorem to the positive cone of $M^d$, we can find a constant, positive-definite matrix $\bar Q$ and $\b^\prime>0$ such that $\L\bar Q=\b^\prime\bar Q$. By the uniqueness of step 1, we have that 
$Q\equiv\bar Q$. 

\noindent{\bf Step 3.} We prove point 2). We saw in lemma 2.1 that $\pc_Q(E,M^d)$ is a convex, compact set of $\mc(E,M^d)$; thus, by Schauder's fixed point theorem, it suffices to show that $\frac{1}{\b}\L^\ast$ brings $\pc_Q(E,M^d)$ into itself. Let 
$\tilde\tau\in\pc_Q(E,M^d)$; we skip the proof that $\frac{1}{\b}\L^\ast\tilde\tau$ is non-negative definite (it follows easily by (2.3) and the definition of the adjoint), but we show that its integral against $Q$ is 1. The first equality below is the definition of the adjoint, the second one is point 1) and the last follows since $\tilde\tau\in\pc_Q(E,M^d)$. 
$$\int_E\left(
Q,\dr\left(
\frac{1}{\b}\L^\ast
\right)   \tilde\tau
\right)_{HS}  =
\int_E\left(
\frac{1}{\b}\L Q,\dr\tilde\tau
\right)_{HS}=
\int_E(Q,\dr\tilde\tau)_{HS}=1  .  $$

\noindent{\bf Step 4.} We prove point 3). Since $\L$ is linear and $B=B^+-B^-$ with 
$B^+,B^-\ge 0$, it suffices to prove (4.17) when $B\in C(E,M^d)$ and $B\ge 0$; in other words, when $B\in\bar C_+$. 

We begin to show that, if $B\in C_+$, then 
$$\theta_+(\L^l B,Q)\tends 0   .   \eqno (4.19)$$
It is clear that (4.19) follows from the three points below. 

\noindent a) For all $\e>0$ there is $a_0>0$ and $\tilde B\in C_+(E,a_0,\nu_0)$ such that 
$\theta^+( B,\tilde B)<\e$. This follows, for instance, since H\"older functions are dense for the $||\cdot||_\infty$ topology. We can also require that $\tilde B\in C_+(E,a,\nu_0)$ for a fixed  $a\ge a_0$. 

\noindent b) If $a$ is large enough, $\theta^{(a,\nu_0)}(\L^l\tilde B,Q)\tends 0$. This follows since $\tilde B\in C_+(E,a,\nu_0)$ and $\L$ is a contraction on $C_+(E,a,\nu_0)$ by step 2. By (1.5), this convergence is exponentially fast. If we apply this argument to 
$\hat B\in C_+(E,a,\nu_0)$ we get the last assertion of the thesis. 

\noindent c) Since $C_+(E,a,\nu_0)\subset C_+$, the definition of hyperbolic distance in section 1 immediately implies that $\theta_+\le\theta^{(a,\nu_0)}$; by the triangle inequality, this implies the first inequality below; the second one comes from the fact that, since 
$\L(C_+)\subset C_+$, then $Lip_{\theta_+}(\L)\le 1$. 
$$\theta_+(\L^l B,Q)\le\theta_+(\L^l B,\L^l\tilde B)+
\theta^{(a,\nu_0)}(\L^l\tilde B,Q)\le$$
$$\theta_+(B,\tilde B)+\theta^{(a,\nu_0)}(\L^l\tilde B,Q)  .  $$
Now the first term on the right is arbitrarily small by point a) and the second one tends to zero by point b). 

We show how (4.19) implies (4.17) when $B\in C_+$. We begin to note that, since 
$\L^\ast\tau=\b\tau$ by point 2) of the thesis, we have for all $l\ge 1$ 
$$\int_E\left(
\left(\frac{1}{\b}\L\right)^l B,\dr\tau
\right)_{HS}  =
\int_E(B,\dr\tau)_{HS}  .  \eqno (4.20)$$
The last formula implies a bound from below on $||\left(\frac{1}{\b}\L\right)^n B||_\infty$; in turn, by lemma 4.1 this implies that the matrices $\left(\frac{1}{\b}\L\right)^n B$ are uniformly positive-definite. Together with (4.19) and lemma 2.2 this implies that there is 
$\a_l>0$ such that 
$$||\a_l\left(\frac{1}{\b}\L\right)^l B-Q||_\infty\tends 0  .  \eqno (4.21)$$
In view of (4.20), this implies that 
$$\a_l\int_E(B,\dr\tau)_{HS}\tends 
\int_E(Q,\dr\tau)_{HS}  .   $$
The right hand side in the formula above is 1 since $\tau\in\pc_Q(G,M^d)$; this implies that $$\a_l\tends\a\colon=\frac{
\int_E(Q,\dr\tau)_{HS}
}{
\int_E(B,\dr\tau)_{HS}
}  .  $$ 
Note that the numerator is 1 since $\tau\in\pc_Q(G,M^d)$; the denominator is different from zero since $B\in C_+$, (2.4) holds and $\tau\in\pc_Q(E,M^d)$. 

Recall that (4.21) and the last formula imply that 
$$\left(
\frac{1}{\b}\L
\right)^l B \tends\frac{1}{\a} Q  $$
uniformly. Now (4.17) follows from the last two formulas. 

The last case is when $B\in\bar C_+\setminus\{ 0 \}$. In this case, we consider $B+\d Id$ for $\d>0$; since $B+\d Id\in C_+$ we have just shown that 
$$\left(\frac{1}{\b}\L\right)^l(B+\d Id)\tends Q\int_G(B+\d Id,\dr\tau)_{HS}  \eqno (4.22)$$
uniformly as $n\tends+\infty$. We saw above that, since $Id\in C_+$, 
$$\left(\frac{1}{\b}\L\right)^l Id\tends Q\int_E(Id,\dr\tau)_{HS}$$
uniformly. Now the thesis follows subtracting the last two formulas. 

\fin

\noindent{\bf Definition.} By point 1) of proposition 4.5, the operator $\L_G$ on $C(G,M^k)$  has a couple eigenvalue-eigenvector which we call $(\b_G,Q_G)$; the operator 
$\L_\Sigma$ on $C(\Sigma,M^k)$ has a couple eigenvalue-eigenvector which we call 
$(\b_\Sigma,Q_\Sigma)$. By point 2) of proposition 4.5 there is a Gibbs measure on $G$, which we call $\tau_G$, and one on $\Sigma$, which we call $\tau_\Sigma$. We shall say that $\kappa_G\colon=(Q_G,\tau_G)_{HS}$ is Kusuoka's measure on $G$ and that 
$\kappa_\Sigma\colon=(Q_\Sigma,\tau_\Sigma)_{HS}$ is Kusuoka's measure on 
$\Sigma$. Since $\tau_G\in\pc_{Q_G}$ and $\tau_\Sigma\in\pc_{Q_\Sigma}$, $\kappa_G$ and $\kappa_\Sigma$ are both probability measures. 

\vskip 1pc

The next lemma shows that there is a natural relationship between these objects. 

\lem{4.6} We have that $\b_G=\b_\Sigma$; we shall call $\b$ their common value. Up to multiplying one of them by a positive constant, we have that 
$Q_\Sigma=Q_G\circ\Phi$. Moreover, $\tau_G=\Phi_\sharp\tau_\Sigma$ and 
$\kappa_G=\Phi_\sharp\kappa_\Sigma$. 

\proof The first equality below comes from the formula after (3.2) and the second one from the definition of $\L_G$ in (3.1). 
$$\L_\Sigma(A\circ\Phi)(x)=
\sum_{i=1}^n \tr D\psi_i|_{\Phi(x)} A_{\psi_i\circ\Phi(x)} D\psi_i|_{\Phi(x)}=
\L_G(A)\circ\Phi (x)
\txt{for all}A\in C(G,M^k)  .   \eqno (4.23)$$
Since $Q_G$ is a fixed point of $\L_G$ on $\frac{C_+(G,a,\nu_0)}{\simeq}$, the formula above implies that $Q_G\circ\Phi$ is a fixed point of $\L_\Sigma$ on 
$\frac{C_+(\Sigma,a,\nu_0)}{\simeq}$. By the uniqueness of proposition 4.5 we get that, up to multiplying one of them by a positive constant,  $Q_\Sigma=Q_G\circ\Phi$. Since 
$$\b_\Sigma Q_\Sigma=\L_\Sigma Q_\Sigma=
\L_\Sigma(Q_G\circ\Phi)=\L_G(Q_G)\circ\Phi=
\b_G Q_G\circ\Phi  $$
we get that $\b_\Sigma=\b_G$. 

We prove the relation between the Gibbs measures. Let $A\in C(G, M^k)$; the first equality below is the definition of the adjoint, the second one follows from the definition of push-forward; the third one comes from (4.23) while the fourth one comes from the fact that 
$\tau_\Sigma$ is an eigenvector of $\L_\Sigma^\ast$; the last one comes again by the definition of push-forward. 
$$\inn{A}{\L_G^\ast(\Phi_\sharp\tau_\Sigma)}=
\inn{\L_G A}{\Phi_\sharp\tau_\Sigma}=
\inn{(\L_G A)\circ\Phi}{\tau_\Sigma}=$$
$$\inn{\L_\Sigma(A\circ\Phi)}{\tau_\Sigma}=
\b\inn{A\circ\Phi}{\tau_\Sigma}=
\b\inn{A}{\Phi_\sharp\tau_\Sigma}  .  $$
Moreover, the fact that $Q_\Sigma=Q_G\circ\Phi$ easily implies that 
$\Phi_\sharp\tau_\Sigma\in\pc_{Q_G}$; since lemma 4.5 implies the uniqueness of the eigenvector of $\L^\ast_G$ in $\pc_{Q_G}$, the last formula implies that 
$\tau_G=\Phi_\sharp\tau_\Sigma$. 

We leave to the reader the easy verification that $\kappa_G=\Phi_\sharp\kappa_\Sigma$. 

\fin

In order to prove that Kusuoka's measure $(Q,\tau)_{HS}$ is ergodic, we need a lemma. 

\lem{4.7} Let $\b>0$ and let $\tau\in\mc^+(E,M^d)$ be as in point 2) of proposition 4.5. Let 
$\fun{A}{E}{M^d}$ be a bounded Borel function. Then, 
$$\int_E (\L A,\dr\tau)_{HS}=
\b\int_E(A,\dr\tau)_{HS}  .  $$

\proof Let us define the measure $t$ as 
$$t\colon=||\tau||+\sum_{i=1}^n(\psi_i)_\sharp(||\tau||)    $$
if we are on $G$; on $\Sigma$ we set 
$$t\colon=||\tau||+\sum_{i=1}^n(a_i)_\sharp(||\tau||)$$
where $\fun{a_i}{(x_0x_1\dots)}{(ix_0x_1\dots)}$. 

By Lusin's theorem there is a sequence $A_k\in C(G,M^d)$ such that $A_k\tends A$ 
$t$-a. e. on $G$; moreover, $||A_k||_\infty$ is bounded. By dominated convergence, this implies that 
$$\int_E(A_k,\dr\tau)_{HS}\tends
\int_E (A,\dr\tau)_{HS}  $$
and
$$\int_E(\L A_k,\dr\tau)_{HS}\tends
\int_E (\L A,\dr\tau)_{HS}   .   $$
Since $A_k$ is continuous, point 2) of proposition 4.5 implies that
$$\int_E(\L A_k,\dr\tau)_{HS}=
\b\int_E(A_k,\dr\tau)_{HS}   .  $$
The thesis follows from the last three formulas. 

\fin

The next lemma recalls some properties of Kusuoka's measure. 

\lem{4.8} Let $E=G$ or $E=\Sigma$; in the first case we set $S=F$, in the second one we set $S=\s$. Then, the following holds. 

\noindent 1) Let $Q$ and $\tau$ be as in proposition 3.2, let $g\in C(E,\R)$ and let 
$A\in C(E,M^k)$. Then we have that 
$$\int_E(g\circ S^l\cdot A,\dr\tau)_{HS}\tends
\int_E(gQ,\dr\tau)_{HS}\cdot
\int_E(A,\dr\tau)_{HS}. \eqno (4.24)$$

\noindent 2) The scalar measures $\kappa_G$ and $\kappa_\Sigma$ defined above are ergodic.

\proof We begin with point 1) on $\Sigma$; we follow [PP]. 

First of all, we note that, if $\fun{h}{\Sigma}{\R}$ is a bounded Borel function, then
$$h\circ \s(ix)=h(x) \txt{for all} i\in(1,\dots,n) \txt{and all} x\in\Sigma  .  \eqno (4.25)$$ 
By (3.2) this implies that, if $A\in C(\Sigma,M^d)$, 
$$[\L_\Sigma(h\circ \s\cdot A)](x)=h(x)(\L_\Sigma A)(x)  \qquad \forall x\in G  .   $$
Integrating against $\tau_\Sigma$, we get the first equality below; the second equality comes from lemma 4.7. 
$$\int_\Sigma\left( 
h\cdot\left(\frac{1}{\b}\L_\Sigma A\right),\dr\tau_\Sigma
\right)_{HS}=
\int_\Sigma\left(
\frac{1}{\b}\L_\Sigma(h\circ \s\cdot A),\dr\tau_\Sigma
\right)_{HS}  =  $$
$$\int_\Sigma(h\circ \s\cdot A,\dr\tau_\Sigma)_{HS}  .  \eqno (4.26)$$
In particular, if $A=Q_\Sigma$ where $Q_\Sigma$ is the eigenfuction of proposition 3.2, we have that 
$$\int_\Sigma h\dr(Q_\Sigma,\tau_\Sigma)_{HS}=
\int_\Sigma h\circ \s\dr(Q_\Sigma,\tau_\Sigma)_{HS}  .  \eqno (4.27)$$
Iterating (4.26) for $h=g\in C(\Sigma,\R)$ we get that 
$$\int_\Sigma\left(
g\left(\frac{1}{\b}\L_\Sigma \right)^kA,\dr\tau_\Sigma
\right)_{HS}=
\int_\Sigma(g\circ \s^k\cdot A,\dr\tau_\Sigma)_{HS}  .  $$
Now (4.24) for $(E,S)=(\Sigma,\s)$ follows from point 3) of proposition 4.5. 

Next, we show (4.24) when $E=G$. Anticipating on lemma 5.3 below, the measure 
$\tau_\Sigma$ is non-atomic. In particular, the countable set $N\subset\Sigma$ on which 
$\Phi$ is not injective is a null set for $\kappa_\Sigma$. Let $g\in C(G,\R)$ and let 
$A\in C(G,M^k)$; the first equality below comes from lemma 4.6; the second one is the definition of push-forward; the third one comes from (1.16), which holds save on a null-set; the limit is (4.24) on $\Sigma$, which we have just proven. The last equality follows again from lemma 4.6.
$$\int_G(g\circ F^l\cdot A,\dr\tau_G)_{HS}=
\int_G(g\circ F^l(x)\cdot A(x),\dr(\Phi_\sharp\tau_\Sigma)(x))_{HS}=$$
$$\int_\Sigma(g\circ F^l\circ\Phi(y)\cdot A\circ\Phi(y), \dr\tau_\Sigma(y))_{HS}=$$
$$\int_\Sigma(g\circ\Phi\circ\sigma^l(y)\cdot A\circ\Phi(y),\dr\tau_\Sigma(y))_{HS}
\tends
\int_\Sigma(g\circ\Phi(y)\cdot Q_\Sigma(y),\dr\tau_\Sigma(y))_{HS}\cdot
\int_\Sigma(A\circ\Phi(y),\dr\tau_\Sigma(y))_{HS}  =  $$
$$\int_G(g(x) Q_G(x),\dr\tau_G(x))_{HS}\cdot
\int_G(A(x),\dr\tau_G(x))_{HS}   .  $$
This is (4.24) for $G$, ending the proof of point 1). 

We prove point 2). First of all, since (4.27) holds for all $h\in C(\Sigma,\R)$ we get that 
$\s_\sharp\kappa_\Sigma=\kappa_\Sigma$, i. e. that $\kappa_\Sigma$ is $\s$-invariant. With the same argument we used for the formula above this implies that 
$$\int_Gh\dr(Q_G,\tau_G)_{HS}=
\int_G h\circ F\dr(Q_G,\tau_G)_{HS}$$
i. e. that $\kappa_G$ is $F$-invariant. 

Now we work on $E$, with $E=G$ or $E=\Sigma$. Setting $A=fQ$ for a continuous function $f$, (4.24) implies that 
$$\int_E g\circ F^l\cdot f\dr (Q,\tau)_{HS}\tends
\int_E g \dr(Q,\tau)_{HS}\cdot
\int_E f \dr(Q,\tau)_{HS}   $$
which implies that $(Q,\tau)_{HS}$ is strongly mixing; in particular, it is ergodic.

\fin

At this stage it is natural to ask whether, when the maps $\psi_i$ are affine, Kusuoka's measure $\kappa_G$ coincides with $(Q_G,\tau_G)_{HS}$; in the remark at the end of section 5 we shall prove that this is the case.

\rm
\vskip 1pc

\vskip 2pc

\centerline{\bf \S 5}
\centerline{\bf The Gibbs property}

\vskip 1pc

In section 1 we defined the cylinder $[x_0\dots x_l]\subset\Sigma$ and the cell  
$[x_0\dots x_l]_G\subset G$.

From now on, we shall suppose that the maps $\psi_i$ are affine and we set 
$$\psi_{x_0\dots x_l}=\psi_{x_0}\circ\dots\circ\psi_{x_l}  .  $$

\vskip 1pc

\noindent{\bf Definition.} Let $\mc^+(G,M^d)$. Following [18], we shall say that $\mu\in\mc^+(G,M^d)$ is a Gibbs measure if there is there are constants $C,D_1>0$ such that, for all $l\ge 1$ and all $x\in G\setminus\tilde N$,
$$e^{-Cl-D_1}\cdot (D\psi_{x_0\dots x_{l-1}})\cdot\mu(G)\cdot\tr (D\psi_{x_0\dots x_{l-1}})\le
\mu ([x_0\dots x_{l-1}]_G)\le $$
$$e^{-Cl+D_1}\cdot (D\psi_{x_0\dots x_{l-1}})\cdot\mu(G)\cdot\tr (D\psi_{x_0\dots x_{l-1}})  .   
\eqno (5.1)_G$$
We  say that $\mu\in\mc^+(\Sigma,M^d)$ is a Gibbs measure if there is there are constants $C,D_1>0$ such that, for all $l\ge 1$ and all $x\in \Sigma$,
$$e^{-Cl-D_1}\cdot (D\psi_{x_0\dots x_{l-1}})\cdot\mu(\Sigma)\cdot\tr (D\psi_{x_0\dots x_{l-1}})\le
\mu ([x_0\dots x_{l-1}])\le $$
$$e^{-Cl+D_1}\cdot (D\psi_{x_0\dots x_{l-1}})\cdot\mu(\Sigma)\cdot\tr (D\psi_{x_0\dots x_{l-1}})  .   
\eqno (5.1)_\Sigma$$

\vskip 1pc

In the formula above, we have not specified at which point we calculate 
$D\psi_{x_0\dots x_{l-1}}$, since $\psi_{x_0\dots x_{l-1}}$ is affine. 

Let $\tau_\Sigma$ be the positive eigenvector of $\L^\ast$ as in lemma 4.5; we briefly prove that $\tau_\Sigma(\Sigma)\not=0$. Let $Q_\Sigma$ be as in proposition 4.5; the inequality below comes from (2.4) and the fact that, for some $\e>0$, 
$\e Q_\Sigma(G)\le Id$ for all $x\in G$ by compactness; the equality comes from the fact that $\tau_\Sigma\in\pc_Q(\Sigma,M^d)$.
$$(Id,\tau_\Sigma(\Sigma))_{HS}\ge
\e\int_\Sigma(Q_\Sigma(x),\dr\tau_\Sigma(x))_{HS}=\e  .  $$

We have the following analogue of proposition 3.2 of [18].

\lem{5.1} Let (F1)-(F4) and (ND) hold; let the maps $\psi_i$ be affine. Let 
$(\b,\tau_\Sigma)$ be as in proposition 4.5. Then for all $l\ge 1$ and all 
$x=(x_0,x_1,\dots)\in\Sigma$ we have 
$$\tau_\Sigma([x_0\dots x_l])=
\frac{1}{\b} (D\psi_{x_0})\cdot
\tau_\Sigma([x_1\dots x_l])\cdot\tr (D\psi_{x_0})   . \eqno (5.2)$$
If $l=0$, (5.2) holds with $\tau_\Sigma(\Sigma)$ instead of $\tau_\Sigma([x_1\dots x_l])$ on the right.

\proof Let 
$x=(x_0,x_1,\dots)\in\Sigma$ be fixed; clearly, we have that 
$$1_{[x_0x_1\dots x_n]}(iz)=
\left\{
\eqalign{
1_{[x_1\dots x_n]}(z)&\txt{if}i=x_0\cr
0&\txt{otherwise.}
}
\right.$$

Let $A\in M^d$ be a fixed, positive semidefinite matrix. The formula above implies the second equality below, while the first one comes  from the fact that $A$ is constant; the third one follows by multiplying and dividing and recalling that 
$1_{[x_0\dots x_l]}(iz)=0$ if $i\not=x_0$; the fourth one from the definition of 
$\L_\Sigma$ in (3.2); the last one follows by lemma 4.7. 
$$(A,\tau_\Sigma([x_1\dots x_l]))_{HS}=
\int_\Sigma(A1_{[x_1\dots x_l]}(z),\dr\tau_\Sigma(z))_{HS}=$$
$$\int_\Sigma\left(
\sum_{i=1}^n A1_{[x_0\dots x_l]}(iz)),\dr\tau_\Sigma(z)
\right)_{HS}=$$
$$\int_\Sigma\Bigg(
\sum_{i=1}^n
\tr D\psi_i\cdot\tr (D\psi_{x_0})^{-1}\cdot 
A1_{[x_0\dots x_n]}(iz)\cdot
(D\psi_{x_0})^{-1}\cdot D\psi_i    ,  
\dr\tau_\Sigma(z)
\Bigg)_{HS}=$$
$$\int_\Sigma\left(
\L_\Sigma\left(
\tr (D\psi_{x_0})^{-1}\cdot A1_{[x_0\dots x_l]}\cdot (D\psi_{x_0})^{-1}
\right)(z),\dr\tau_\Sigma(z)
\right)_{HS}   =$$
$$\b\int_{[x_0\dots x_l]}\left(
\tr (D\psi_{x_0})^{-1}\cdot A\cdot (D\psi_{x_0})^{-1},\dr\tau_\Sigma
\right)_{HS}  .   
\eqno (5.3)$$
After transposition, (5.3) implies that
$$\b
(A, D\psi_{x_0}^{-1}\cdot\tau_\Sigma([x_0\dots x_l])\cdot\tr D\psi_{x_0}^{-1})_{HS} =
(A,\tau_\Sigma([x_1,\dots x_l]))_{HS}  .  $$
Letting $A$ vary among the one-dimensional projections we get that
$$\b\cdot D\psi_{x_0}^{-1}\cdot\tau_\Sigma([x_0\dots x_l])\cdot\tr D\psi_{x_0}^{-1}=
\tau_\Sigma[x_1\dots x_l]  .  $$
To get (5.2) it suffices to multiply the formula above by 
$$\frac{1}{\b}\cdot (D\psi_{x_0}) \txt{on the left and by} 
\tr(D\psi_{x_0})\txt{on the right.}$$

\fin

\cor{5.2} Let (F1)-(F4) and (ND) hold; let us suppose that the maps $\psi_i$ are affine and  let $(\b,\tau_\Sigma)$ be as in proposition 4.5. then, $\tau_\Sigma$ is a Gibbs measure for the constant $C=\log\b$. 

\proof Iterating the right hand side of (5.2) and using the chain rule we get the following.
$$\tau_\Sigma([x_0\dots x_l])= \frac{1}{\b}\cdot
(D\psi_{x_0})\tau_\Sigma([x_1\dots x_l]) \tr(D\psi_{x_0})=$$
$$\frac{1}{\b^2}\cdot (D\psi_{x_0x_1})\tau_\Sigma([x_2\dots x_l]) \tr(D\psi_{x_0x_1})=$$
$$\dots=$$
$$\frac{1}{\b^{l}}\cdot (D\psi_{x_0x_1\dots x_{l-1}})\tau_\Sigma([x_l]) 
\tr(D\psi_{x_0x_1\dots x_{l-1}})   =  $$
$$\frac{1}{\b^{l+1}}\cdot (D\psi_{x_0\dots x_l})\tau_\Sigma(\Sigma) \tr(D\psi_{x_0\dots x_l})   .  $$

\fin

\lem{5.3} Let the maps $\psi_i$ satisfy (F1)-(F4) and let (ND) hold. Then, we have the following. 

\noindent 1) The measure $\tau_\Sigma$ is positive on open sets. 

\noindent 2) The measures $\tau_\Sigma$ and $\tau_G$ are non-atomic. 

\proof We begin with point 1) for $\tau_\Sigma$. It suffices to show that, for all cylinders 
$[x_0\dots x_l]\subset\Sigma$, the matrix $\tau_\Sigma[x_0\dots x_l]$ is not zero. We get from (5.3) that 
$$(Id,\tau_\Sigma[x_1\dots x_l])_{HS}=
\b\int_{[x_0\dots x_l]}(
\tr(D\psi_{x_0}|_{\Phi(x)})^{-1}\cdot (D\psi_{x_0}|_{\Phi(x)})^{-1},\dr\tau_\Sigma(x)
)_{HS}   .  $$
This easily implies that, if $\tau_\Sigma[x_1\dots x_l]$ is not zero, then also 
$\tau_\Sigma[x_0\dots x_l]$ is not zero. Iterating, we see that $\tau_\Sigma[x_0\dots x_l]$ is not zero if $\tau_\Sigma(\Sigma)$ is not zero, a fact we showed before stating this lemma. 

As for point 2), we begin to recall the standard proof that $\kappa_\Sigma$ is non-atomic. By point 2) of lemma 4.8, $\kappa_\Sigma$ is ergodic; let us suppose by contradiction that it has an atom $\{ \bar x \}$. We are going to show that $\kappa_\Sigma(\{ \bar x \})=1$ and, consequently, $\kappa_\Sigma(\{ \bar x \}^c)=0$. This will be the contradiction, since by point 1) $\tau_\Sigma$  is positive on open sets. 

First of all, let us suppose that $\bar x$ is a periodic orbit of period $q$ and let us set 
$$A=\bigcup_{l\ge 0}\s^{-lq}(\{ \bar x \})   .   $$
Clearly, $A$ is $\s^q$-invariant, i. e. 
$\s^{-q}(A)\subset A$. Since $\s$ preserves $\kappa_\Sigma$, we have that 
$\kappa_\Sigma(\s^{-lq}\{ \bar x \})=\kappa_\Sigma(\{ \bar x \})$; since $\s^q$ fixes 
$\bar x$, we see that $\bar x \in\s^{-q}(\{ \bar x \})$. This implies that  
$\tau_\Sigma$ on $\s^{-q}(\{ \bar x \})$ concentrates on $\{ \bar x \}$; iterating, we get that 
$\tau_\Sigma$ on $\s^{-lq}(\{ \bar x \})$ concentrates on $\{ \bar x \}$. This and the definition of $A$ easily imply that $\kappa_\Sigma(A)=\kappa_\Sigma(\{ \bar x \})$ and that 
$\kappa_\Sigma(A\setminus\s^{-q}(A))=0$. By ergodicity, this implies that $\kappa_\Sigma(\{ \bar x \})=1$. 

The second case is when $\bar x$ has an antiperiod, say of length $l$. We consider 
$\tilde x=\s^l(\bar x)$, which is periodic. Since $\bar x\in\s^{-l}(\tilde x)$, invariance implies that $\kappa_\Sigma(\tilde x)>0$; now the same argument as above applies.  

The last case is when $\bar x$ is not periodic; then, it is easy to see that the sets 
$\s^{-l}(\{ \bar x \})$ are all disjoint. Since they have the same measure, we get that 
$\kappa_\Sigma(\{ \bar x \})=0$, i. e. that $\{ \bar x \}$ is not an atom. 

In order to show that $\kappa_G$ is non-atomic, it suffices to recall three facts: that 
$\tau_\Sigma$ is non-atomic, that $\tau_G=\Phi_\sharp\tau_\Sigma$ by lemma 4.6 and that $\Phi$ is finite-to-one. 

\fin

\noindent{\bf End of the proof of theorem 1.} Points 1) and 2) come from proposition 4.5. The self-similarity of point 4) comes from lemma 3.1. The ergodicity of point 3) is point 2) of lemma 4.8; mutual absolute continuity in one direction follows from (2.8), in the other one is trivial.  For point 5) we begin to note that, by the definition of $\Phi$, 
$[x_0\dots x_l]\subset\Phi^{-1}([x_0\dots x_l]_G)$; the points of 
$\Phi^{-1}([x_0\dots x_l]_G\setminus [x_0\dots x_l]$ are those with multiple codings, which we have seen in section 1 to be a countable set. Since $\tau_\Sigma$ is non-atomic, we get that 
$$\tau_G(
\Phi^{-1}([x_0\dots x_l]_G\setminus [x_0\dots x_l]
)=0  .   $$
Since $\tau_G=\Phi_\sharp\tau_\Sigma$ by lemma 3.4, the last formula implies that 
$$\tau_\Sigma([x_0\dots x_l])=\tau_G([x_0\dots x_l]_G)   .  $$
Since $\tau_\Sigma$ has the Gibbs property by corollary 5.2, we are done. 

\fin

\noindent{\bf Remark.} In corollary (4.2) we have supposed that the maps $D\psi_i$ are constant, which is the case of Kusuoka's paper [12].  We prove that, up to multiplication by a positive constant, $(Q_G,\tau_G)_{HS}$ coincides with Kusuoka's measure $\kappa$. 

When $D\psi_i$ is constant, point 1) of proposition 4.5 implies that $Q_G$ is constant too and solves 
$$Q_G=\frac{1}{\b}\sum_{i=1}^n\tr D\psi_i Q  D\psi_i  .   \eqno (5.4)$$
By point 2) of lemma 3.1 we have that 
$$\b\ec_\tau(f,g)=
\sum_{i=1}^n\ec_\tau(f\circ\psi_i,g\circ\psi_i)=$$
$$\sum_{i=1}^n\int_G(
\tr D\psi_i\cdot\nabla f|_{\psi_i(x)},\dr\tau\tr D\psi_i\cdot\nabla g|_{\psi_i(x)} 
)   . $$
If we choose as $f$ and $g$ two linear functions and we recall that $D\psi_i$ is a constant matrix, we see that the last formula implies that
$$\b\tau_G(G)=\sum_{i=1}^n
D\psi_i\cdot\tau_G(G)\cdot\tr D\psi_i  .   \eqno (5.5)$$
Kusuoka's measure $\kappa$ is defined by the following formula: if $x=(x_0x_1\dots)$, then 
$$\kappa([x_0\dots x_l]_G)=
\frac{1}{\b^l}(
Q,\tr (D\psi_{x_0\dots x_l}) \hat Q  \tr(D\psi_{x_0\dots x_l})
)_{HS}$$
where $Q$ solves (5.4) and $\hat Q$ solves (5.5). Since the solution to both equations is unique by Perron-Frobenius (up to multiplication by a constant, of course), the last formula and corollary (5.2) imply that $\kappa_G\colon=(Q_G,\tau_G)$ coincide with $\kappa$.

\vskip 2pc

\vskip 2pc
\centerline{\bf References}

\noindent [1] M. T. Barlow, R. F. Bass, The construction of Brownian motion on the Sierpinski carpet, Ann. IHP,  {\bf 25}, 225-257, 1989.

\noindent [2] M. T. Barlow, E. A. Perkins, Brownian motion on the Sierpiski gasket, Probab. Th. Rel. Fields, {\bf 79}, 543-623, 1988. 

\noindent [3] R. Bell, C. W. Ho, R. S. Strichartz, Energy measures of harmonic functions on the Sierpinski gasket, Indiana Univ. Math. J. {\bf 63}, 831-868, 2014. 

\noindent [4] G. Birkhoff, Lattice theory, Third edition, AMS Colloquium Publ., Vol. XXV, AMS, Providence, R. I., 1967.

\noindent [5] M. Fukushima, Y. Oshima, M. Takeda, Dirichlet forms and symmetric Markov processes, De Gruyter, G\"ottingen, 2011. 

\noindent [6] S. Goldstein, Random walks and diffusions on fractals, in H. Kesten (ed), Percolation theory and ergodic theory of infinite particle systems, IMA vol. Math. Appl., 
{\bf 8}, Springer, New York, 121-129, 1987.

\noindent [7] A. Johansson, A. \"Oberg, M. Pollicott, Ergodic theory of Kusuoka's measures, J. Fractal Geom., {\bf 4}, 185-214, 2017.  

\noindent [8] N. Kajino, Analysis and geometry of the measurable Riemannian structure on the Sierpinski gasket, Contemporary Math. {\bf 600}, Amer. Math. Soc., Providence, RI, 2013.  

\noindent [9] J. Kigami, Analysis on fractals, Cambridge tracts in Math., {\bf 143}, Cambridge Univ. Press, Cambridge, 2001. 

\noindent [10] P. Koskela, Y. Zhou, Geometry and Analysis of Dirichlet forms, Adv. Math., 
{\bf 231}, 2755-2801, 2012. 

\noindent [11] S. Kusuoka, A diffusion process on a fractal, in: K. Ito and N. Ikeda (eds.), Probabilistic methods in Mathematical Physics, Academic Press, Boston, MA, 251-274, 1987. 

\noindent [12] S. Kusuoka, Dirichlet forms on fractals and products of random matrices, Publ. Res. Inst. Math. Sci., {\bf 25}, 659-680, 1989.

\noindent [13] R. Ma\~n\'e, Ergodic theory and differentiable dynamics, Berlin, 1983.

\noindent [14] U. Mosco, Composite media and asymptotic Dirichlet forms, J. Functional Analysis, {\bf 123}, 368-421, 1994.

\noindent [15] U. Mosco, Variational fractals, Ann. Scuola Norm. Sup. Pisa Cl. Sci. (4), 
{\bf 25}, 683-712, 1997.

\noindent [16] R. Peirone, Existence of self-similar energies on finitely ramified fractals, J. Anal. Math., {\bf 123}, 35-94, 2014.

\noindent [17] R. Peirone, Convergence of Dirichlet forms on fractals, mimeographed notes. 

\noindent [18] W. Perry, M. Pollicott, Zeta functions and the periodic orbit structure of hyperbolic dynamics, Asterisque, {\bf 187-188}, 1990. 



\noindent [19] M. Viana, Stochastic analysis of deterministic systems, mimeographed notes. 

\noindent [20] J. Claude Yoccoz, Hyperbolic dynamics, mimeographed notes, 1991.

\end

\centerline{\bf \S 6}
\centerline{\bf The cone-field}

\vskip 1pc

We begin by recalling the definition and some properties of the cone fields of Dynamical Systems; we shall use them in lemma 6.3 below to show that the measure $\tau$ of proposition 4.5 has the form 
$$\tau=P_{z_x}||\tau||$$
where $z_x$ is a Borel field of vectors and $P_{z_x}$ is the orthogonal projection on 
$z_x$. We are interested in fractals in which $\tau$ relates to a one-dimensional field of projections which, essentially, are the projections on the "less expanding" directions of $F$; the notion can be made precise using $\s-\s^\prime$-hyperbolicity as in [20]. 

\vskip 1pc

\noindent{\bf Definitions.} We take from [20] the classical definition of cone-field with a one dimensional core. For $v\in\R^d\setminus\{ 0 \}$, we denote by $P_v$ and by $P_v^\perp$ the orthogonal projections on $v$ and $v^\perp$ respectively: 
$$P_v w=\frac{(v,w)}{||v||^2}v,\qquad P_v^\perp w=w-P_vw  .  $$
If $v\in\R^d\setminus\{ 0 \}$ and $\l>0$, we define 
$$Cone(v,\l)=
\{
w\in\R^d\st ||P_v^\perp w||\le\l ||P_v w||
\}   \eqno (6.1)$$
and 
$$Cone^+(v,\l)=\{
w\in Cone(v,\l)\st (v,w)\ge 0
\}   .  \eqno (6.2)$$
Note that $Cone^+(v,\l)$ is one of the convex cones we considered in section 1. 

Let $H\subset G$ be a Borel set such that $F(H)=H$; let $\fun{v}{G}{\R^d\setminus\{ 0 \}}$ be a Borel vector field and let $\fun{\l}{G}{(0,+\infty)}$ be Borel; we associate to each 
$x\in H$ the cone $Cone(v_x,\l_x)$. If there is $D_1>0$ such that 
$$\frac{1}{D_1}\le \l_x\le D_1    \eqno (6.3)$$
we say that $Cone(v_x,\l_x)$ is a bounded cone-field over $H$.  

We recall a standard proposition (see for instance section 4.7 of [20]).

\lem{6.1} Let $H\subset G$ be an invariant Borel set and let $Cone(v_x,\l_x)$ be a bounded cone-field over $H$. We suppose that there is $b\in (0,1)$ such that 
$$(DF)^{-1}(x)Cone^+(v_{F(x)},\l_{F(x)})\subset
Cone^+(v_x,b\l_x)   \txt{for all} x\in H    .    \eqno (6.4)$$

Then, the following holds.

\noindent 1) There is a Borel vector field $\fun{z}{H}{\R^d}$ and a Borel function 
$\fun{\g}{H}{\R}$ such that, for all $x\in H$, 
$$(x,z_x)\   \subset Cone(v_x,\l_x)   $$
and 
$$(DF)^{-1}(x)z_{F(x)}=\g_x z_x  .  $$

\noindent 2) Moreover, we have that 
$$\frac{(DF^l)^{-1}|_{F^l(x)}w}{||(DF^l)^{-1}|_{F^l(x)}w||}\tends 
P_{z_x}\frac{w}{||w||}   \eqno (6.5)$$
uniformly for $(x,w)\in H\times B(0,1)$. 

\rm

\vskip 1pc

\noindent{\bf Sketch of the proof.} We only sketch the standard proof. On each 
$Cone^+(v_x,\l_x)$ we put the hyperbolic distance $\theta_x$ and we consider the space 
$$C=\prod_{x\in H}Cone^+(v_x,\l_x)$$
with the product topology. By Tychonov's theorem,
$$C_b\colon=\prod_{x\in H}Cone^+(v_x,b\l_x)  $$
is compact in $C$. Since the map 
$$\fun{F_\ast}{B(H,\R^d)}{B(H,\R^d)}$$
$$(F_\ast v)(x)=(DF)^{-1}(x)v_{F(x)}      $$
lands in $C_b$, Schauder's fixed point theorem implies that there is $z\in C$ satisfying point 1). 

By (6.4), (1.5) and the fact that the cone-field is bounded, we can see that 
$${\rm diam}_{\theta_x}(
(DF^l)^{-1}(x)Cone(v_{F^l(x)},\l_{F^l(x)})
)  \tends 0  .  $$
It is easy to see that this implies point 2). 

\fin

\lem{6.2} Let us suppose that $G$ satisfies (F1)-(F4) and (ND); let the maps $\psi_i$ be affine. Let us also suppose that there is an invariant Borel set $H$ of full $||\tau||$-measure and a bounded cone-field over $H$ satisfying hypotheses 1) and 2) of lemma 6.1; let 
$\fun{z}{G}{\R^d}$ be the vector field given by this lemma and let $\tau$ be the Gibbs measure of proposition 4.5. Then, 
$$\tau=P_{z_x}||\tau||  .  $$

\proof We begin with point 1). Let $\tau=T_x||\tau||$ be the polar decomposition of $\tau$. Recall that the total variation is defined with respect to $||\cdot||_{HS}$ on $M^d$; as a consequence, $||T_x||_{HS}=1$ for $||\tau||$-a. e. $x\in G$. By Rokhlin's theorem (see for instance [19] or [20], or theorem 0.5.1 of [13]) we have that, for $||\tau||$-a. e. $x\in G$, the equality below holds. 
$$T_x=\lim_{l\tends+\infty}\frac{
\tau(G([x_0\dots x_l]))
}{
||\tau||(G([x_0\dots x_l]))
}   .  $$
Since $||T_x||_{HS}=1$, the last formula implies that 
$$\lim_{l\tends+\infty}\frac{
||\tau||(G([x_0\dots x_l]))
}{
||\tau(G([x_0\dots x_l]))||_{HS}
}  =1     $$
for $||\tau||$-a. e. $x\in G$. From the last two formulas we deduce that 
$$T_x=\lim_{l\tends+\infty}\frac{
\tau(G([x_0\dots x_l]))
}{
||\tau(G([x_0\dots x_l]))||_{HS}
}    .  $$
By the formula of corollary 5.2 we get that 
$$\frac{
\tau(G([x_0\dots x_l]))
}{
||\tau(G([x_0\dots x_l]))||_{HS}
}  =
\frac{\tr (DF^l)^{-1}(x)\tau(G) (DF^l)^{-1}(x)}{||\tr (DF^l)^{-1}(x)\tau(G) (DF^l)^{-1}(x)||_{HS}}   .  $$
Now the thesis comes from point 2) of lemma 6.1.

\fin

\noindent{\bf Proof of theorem 1.} Points 1), 2) and 3) of theorem 1 come from proposition 4.5. Point 4) is lemma 3.1. Point 5) follows from corollary 5.2, point 6) from lemma 6.2.  

\fin